\documentclass[11pt]{article}

\usepackage[margin=1in]{geometry}
\usepackage{amsmath,amsfonts,amsthm,amssymb}
\usepackage{subcaption}
\usepackage{graphicx}
\graphicspath{{figures/}}
\usepackage{color}
\usepackage{url}
\usepackage{soul}

\usepackage{verbatim}%for large commenting
\usepackage{hyperref}
\title{Regularity of Solutions for the Nonlocal Wave Equation on Periodic Distributions}
\author{Thinh Dang, Bacim Alali, and Nathan Albin\\
\
\footnotesize{Department of Mathematics, Kansas State University, Manhattan, KS}}
\date{}

\newcommand{\Rn}{\mathbb{R}^n}

\newcommand{\Ldel}{L^{\delta,\beta}}

\newcommand{\cdel}{c^{\delta,\beta}}
\newcommand{\mdel}{m^{\delta,\beta}}

\newcommand{\T}{\mathbf{T}}

\newtheorem{theorem}{Theorem}
\newtheorem{definition}{Definition}
\newtheorem{lemma}{Lemma\textbf{}}
\newtheorem{proposition}{Proposition\textbf{}}
\newtheorem{cor}{Corollary}

\theoremstyle{remark}
\newtheorem{remark}{Remark}

\begin{document}
 \title{Regularity of Solutions for the Nonlocal Wave Equation on Periodic Distributions
 %\thanks{Received date, and accepted date (The correct dates will be entered by the editor).}
 }

          %For each author, make a block with the following macros:

          \author{Thinh Dang	\thanks{
          tdang@ksu.edu, Department of Mathematics, Kansas State University, Manhattan, KS.}
          \and Bacim Alali \thanks{bacimalali@math.ksu.edu, Department of Mathematics, Kansas State University, Manhattan, KS.}
          \and  Nathan Albin \thanks{albin@ksu.edu, Department of Mathematics, Kansas State University, Manhattan, KS.}
          \thanks{This project is based upon work supported by the National Science Foundation under Grant No. 2108588.}}
%\footnote{Department of Mathematics, Kansas State University, Manhattan, KS.}

         %\pagestyle{myheadings} \markboth{Regularity of Solutions for the Nonlocal Wave Equation on Periodic Distributions}{T. DANG, B. ALALI \and N. ALBIN} 
         \maketitle

          \begin{abstract}
This work addresses the regularity of solutions for a nonlocal wave equation over the space of periodic distributions. The spatial operator for the nonlocal wave equation is given by a nonlocal Laplace operator with a compactly supported integral kernel. We follow a unified approach based on the Fourier multipliers of the nonlocal Laplace operator, which allows the study of regular as well as distributional solutions of the nonlocal wave equation, integrable as well as singular kernels, in any spatial dimension. In addition, the results extend beyond operators with singular kernels to nonlocal-pseudo differential operators. We present results on the spatial and temporal regularity of solutions in terms of regularity of the initial data or the forcing term. Moreover, solutions of the nonlocal wave equation are shown to converge to the solution of the classical wave equation for two types of limits: as the spatial nonlocality vanishes or as the singularity of the integral kernel approaches a certain critical singularity that depends on the spatial dimension.
          \end{abstract}
% \begin{keywords}  nonlocal Laplacian; wave equation; distributions; Fourier multipliers; peridynamics.
% \end{keywords}

% \begin{AMS} 45A05; 45P05; 47G10; 45M15.
% \end{AMS}
\section{Introduction}
This work is focused on the regularity of  solutions to the nonlocal wave equation
\[
\begin{cases}
    u_{tt}(x,t)=L_\gamma u(x,t)+b(x),&\quad x\in \T^n,\; t>0\\
    u(x,0)=f(x),& \quad x\in\T^n\\
    u_t(x,0)=g(x),& \quad x\in\T^n
\end{cases},
\]
%where the divergence operator in the local counterpart is replaced by 
where ${\T}^n$ is the periodic torus in $\mathbb{R}^n$. Here
$L_\gamma$ is a nonlocal Laplace operator of the form
\begin{align}
\label{eq:Lgamma}
    L_{\gamma}u(x)=\int(u(y)-u(x))\gamma(y,x)dy,
\end{align}
with $\gamma$ being a symmetric kernel with compact support.
%, i.e. $\gamma(x,y)=\gamma(y,x)$ for any $x,y\in\mathbb{R}^n$. 
This type of operator, in the case when $u$ is a vector field and $\gamma$ being a second order tensor field, originated in  peridynamics \cite{Silling2000,silling2007peridynamic,silling2010linearized}. For the case of scalar fields $u$, this operator has been introduced in Nonlocal Vector Calculus \cite{nonlocal_calc_2013}. Regularity of solutions for nonlocal equations associated with operators of form \eqref{eq:Lgamma} has been studied under different settings. The work in \cite{FossRadu} studies a Dirichlet-type (volume-constraint)   problem for scalar nonlocal equations and the work in \cite{alali2021fourier} studies the  periodic nonlocal Poisson equation. Peridynamics Dirichlet-type constraint problems have been studied in the works  \cite{kassmann2019solvability,MengeshaDu,FossRaduYu}. The work in \cite{Beyer} addresses nonlocal wave equations on bounded domains and the work in \cite{mengesha2020solvability} discusses the solvability of peridynamics equations in  $\Rn$. Regularity of solutions for peridynamics equations, associated with non-compactly supported integral kernels, is addressed in \cite{Alimov2014}. There are several numerical and computational studies of nonlocal wave equations and peridynamics equations  including the works in \cite{guan2015stability, du2018numerical, coclite2020numerical, du2019nonlocal,d2020numerical,alali2020fourier}.
%,d2020numerical,du2016asymptotically,jafarzadeh2020efficient,du2017fast,du2018numerical}.

%This type of operator has been  studied in the context of equilibrium equations or scalar Dirichlet peridynamic equations on bounded domains, see for example ~\cite{Beyer,FossRadu,FossRaduYu,kassmann2019solvability,mengesha2020solvability,MengeshaDu} . In the context of nonlocal wave equations, with second-rank kernel $\gamma$ and tensor values paired with vector fields, the corresponding operators are defined to be peridynamic operator~\cite{Silling2000,silling2010linearized}. The nonlocal wave equation with linear peridynamics were studied in~\cite{Alimov2014} where the equation was converted into an operator valued Volterra integral equation of second kind with the differential-integral operator kernel. The solution regularity was then derived through the well-posedness of the Volterra integral equation. When the kernel $\gamma$ and field $u$ are scalar-valued, the corresponding nonlocal wave equation was studied in both bounded settings~\cite{Kun2010,Quing2015} and unbounded settings~\cite{QiangDu2018}.   

In this work, we focus on integral kernels  of the form
\begin{align*}
    \gamma(x,y)=\frac{\cdel}{\|y-x\|^\beta}\chi_{B_\delta(x)}(y)
\end{align*}
where $\chi_{B_\delta(x)}$ denotes the indicator function of the ball centered at $x$ with radius $\delta>0$, $\cdel$ a scaling constant, and the kernel exponent $\beta<n+2$. 

For $\delta>0$ and $\beta< n+2$, the nonlocal Laplace operator $\Ldel$ is defined by
\begin{equation}\label{eq:nonlocal_laplacian}
  \Ldel u(x) = \cdel\int_{B_\delta(x)}\frac{u(y)-u(x)}{\|y-x\|^\beta}\;dy=\cdel\int_{B_\delta(0)}\frac{u(x+z)-u(x)}{\|z\|^\beta}\;dz, 
\end{equation}
where the scaling constant $\cdel$ is 
%chosen in such a way that for the function $u(z):=\|z\|^2$, \\$\Ldel u = 2 n = \Delta u$, or equivalently, it is 
given by
\begin{eqnarray*}
\label{eq:cdel-integral}
\cdel 
\label{eq:cdel-explicit}
=\frac{2(n+2-\beta)\Gamma\left(\frac{n}{2}+1\right)}
    {\pi^{n/2}\delta^{n+2-\beta}}.
\end{eqnarray*}

In~\cite{alali2021fourier}, the Fourier multiplier  of this operator is shown to have a hypergeometric representation, which then is used to extend the multiplier, and consequently the operator, to  the case in which $\beta \ge n+2$, with $\beta\ne n+4,n+6,\ldots$ %The operator itself was also defined in this case as a pseudo operator to match the formula of the multiplier. 
We exploit the Fourier multipliers results that were developed in~\cite{alali2021fourier} to study the regularity of the  nonlocal wave equation over the space of periodic distributions $H^s(\T^n)$, with $s\in\mathbb{R}$. 

The organization  and the main contributions of this article are described  as follows. A review of some results on the multipliers $\mdel$ of the operator $\Ldel$ are provided in Section \ref{sec:overview}. The operator $\Ldel$, with $\beta\ge n+2$, which has been defined as a pseudo operator in \cite{alali2021fourier}, is shown to have an explicit mixed integral-differential representation for the case when $n+2\le\beta<n+4$ in Section \ref{sec:L-explicit}. In addition, this formula is  shown to be consistent with the original integral formula for $\Ldel$, given by \eqref{eq:nonlocal_laplacian}, in the case when $\beta<n+2$. In Section \ref{sec:monotonicity}, we prove a result on the monotonicity of the multipliers $\mdel$ as a function of the kernel exponent $\beta$, which will be used in Sections \ref{sec:convergence1} and \ref{sec:convergence2} to show convergence of solutions results in the limit as $\beta\rightarrow n+2$. We consider the periodic nonlocal wave equation with initial data, in Section \ref{sec:Homogeneous}, or with a forcing term, in Section \ref{sec:forcing}. We prove spatial regularity results for these equations over the space of periodic distributions. Temporal regularity of solutions of the nonlocal wave equation is discussed in Sections \ref{sec:temporal1} and \ref{sec:temporal2}. The time derivative is defined as a Gateaux derivative when the solutions are distributions and we use the classical time derivative when the solutions are in $L^2(\T^n)$ with respect to the spatial variable. In Sections \ref{sec:convergence1} and \ref{sec:convergence2}, convergence of solutions of the nonlocal wave equation to the solution of the classical wave equation are proved for two types of limiting behavior: 
as $\delta \to 0^+$, with $\beta< n+4$ is being fixed,  or as $\beta\to n+2$, with $\delta>0$ is being fixed.
%\vspace{5cm}
%The nonlocal Laplace operator is defined in Section 2.1 along with its corresponding multipliers for both case $\beta<n+2 $ and $n+2\le \beta<n+4$.  The remaining of section 2 mentions some Lemmas, Propositions and Theorems that would be used throughout the article. Section 3 provides us with an explicit representation for the nonlocal Laplace operator for the general case $\beta<n+4$. This formula is also shown to be consistent with the formula for the case $\beta<n+2$ that was previously well-known. It is also expected from the multiplier that the case $\beta=n+2$, the operator coincides perfectly with the famous classical Laplace operator. The monotonicity of the multipliers $\mdel(\nu)$ as a function of the kernel exponenet $\beta$ is discussed in section 3.2. Section 3.3 mentions the periodic set up of the problem with the eigenvalues of the operator in the domain. In section 3.4, we define the Sobolev space for periodic distributions, the temporal derivative and the result when applying the nonlocal operator to those distributions.   Then we dedicate Section 4 and 5 looking at the regularity spatially and temporally of the solution to the nonlocal Laplace wave equation in 2 separate settings: homogeneous with nonzero initial conditions and inhomogeneous with external force term and zero initial conditions. We also discuss how the solutions converge to the solution of the classical cases under special limitations. Finally, in Section 6, we combine the previous two results for the most general case of nonlocal wave equation.    

\section{Overview and some results on the nonlocal Laplacian}
\label{sec:overview}
In this section, we present some previous as well as new results on the nonlocal Laplace operator \eqref{eq:nonlocal_laplacian}. In addition, we provide some results on the space of periodic distributions that will be used in subsequent sections.

%e following subsections, we provide an overview on the nonlocal Laplacian and its Fourier multipliers.
\subsection{The nonlocal Laplace operator}
\label{sec:prelim}
\begin{comment}
For $\delta>0$ and $\beta< n+2$, the nonlocal Laplace operator $\Ldel$ is defined as
\begin{equation}\label{eq:nonlocal_laplacian}
  \Ldel u(x) = \cdel\int_{B_\delta(x)}\frac{u(y)-u(x)}{\|y-x\|^\beta}\;dy=\cdel\int_{B_\delta(0)}\frac{u(x+z)-u(x)}{\|z\|^\beta}\;dz, 
\end{equation}
where the scaling constant $\cdel$ is 
%chosen in such a way that for the function $u(z):=\|z\|^2$, \\$\Ldel u = 2 n = \Delta u$, or equivalently, it is 
given by
\begin{eqnarray*}
\label{eq:cdel-integral}
\cdel 
\label{eq:cdel-explicit}
=\frac{2(n+2-\beta)\Gamma\left(\frac{n}{2}+1\right)}
    {\pi^{n/2}\delta^{n+2-\beta}}.
\end{eqnarray*}

\noindent
\end{comment}
The Fourier multiplier $\mdel$ of the operator $\Ldel$, in \eqref{eq:nonlocal_laplacian}, are defined through  the Fourier transform by
\begin{equation*}
\label{eq:multiplier-def}
   \widehat{\Ldel u} = \mdel \; \hat{u}.
\end{equation*}
where  $\mdel$ is of the form
\begin{equation}\label{eq:multiplier-cosine0}
	\mdel(\nu) = \cdel\int_{B_\delta(0)}\frac{\cos(\nu\cdot z)-1}{\|z\|^\beta}\;dz.
\end{equation}

\noindent
In subsequent sections, we will make use of the following results, which are proved in \cite{alali2021fourier}.
\begin{theorem}\label{thm:multipliers-2F3}
Let $n\ge 1$, $\delta>0$ and $\beta<n+2$.  The Fourier multipliers $\mdel$ have the following representation
\begin{equation}\label{eq:multiplier-general}
	\mdel(\nu) = -\|\nu\|^2\,_2F_3\left(1,\frac{n+2-\beta}{2};2,\frac{n+2}{2},\frac{n+4-\beta}{2};-\frac{1}{4}\|\nu\|^2\delta^2\right).
\end{equation}
\end{theorem}
\noindent The hypergeometric function $_2F_3$ on the right hand side is well-defined for any $\beta\ne n+4,n+6,\cdots$, hence, using \eqref{eq:multiplier-general}, the definition of the multipliers is extended to the case when $\beta\ge n+2$ with $\beta\ne n+4, n+6,\cdots$. 
Consequently, the operator $\Ldel$ is extended to these larger  values of $\beta$ using the Fourier transform.
%define the Fourier multipliers using the same formula for the case where 
%Then, in these cases, the Nonlocal Laplace operator is extended by being defined through its multipliers as a pseudo-differential operator. Specifically, 
In particular, for the case  when $\beta=n+2$, and since $m^{\delta,n+2}(\nu)$ is equal to $-\|\nu\|^2$, the extended operator $\Ldel$ coincides with the classical Laplace operator $\Delta$.

For clarity, we reiterate the fact  that the operator $\Ldel$ is an integral operator for $\beta<n+2$, having an integral kernel when $\beta<n$ and a singular kernel when $n\le \beta<n+2$. When $\beta>n+2$, with $\beta\ne n+4, n+6, \ldots$, the operator $\Ldel$ can be recognized as an integral-differential operator, as we show in Section~\ref{sec:L-explicit}. In particular, when $n+2<\beta<n+4$, the operator $\Ldel$ can be interpreted as a  super-diffusion operator in the setting of heat diffusion or conduction, as discussed in \cite{alali2020fourier}.
\begin{cor}\label{convergenge_L}
Let $n\ge 1$ and $\nu\in\Rn$.  Then the  Fourier multipliers $\mdel$   converge to the Fourier multipliers of the Laplacian  as follows
\begin{equation}\label{col:delta}
    \lim_{\delta\rightarrow 0^+}\mdel(\nu)=-\|\nu\|^2, \;\;\mbox{ for } \;\; \beta\in\mathbb{R}\setminus\{n+4,n+6,n+8,\ldots\},
\end{equation}
and
\begin{equation}\label{col:beta}
    \lim_{\beta\rightarrow n+2} \mdel(\nu)=-\|\nu\|^2, \;\;\mbox{ for }\;\; \delta>0.
\end{equation}
\end{cor}
The asymptotic behavior of the multipliers for large $\nu$ is described through the following result.
\begin{theorem}\label{thm:asymptotics}
Let $n\ge 1$, $\delta>0$ and $\beta\in\mathbb{R}\setminus\{n+2,n+4,n+6,\ldots\}$.  Then, as $\|\nu\|\to\infty$,
\begin{equation*}
\mdel(\nu) \sim 
\begin{cases}
-\frac{2n(n+2-\beta)}{\delta^2(n-\beta)}
+ 2\left(\frac{2}{\delta}\right)^{n+2-\beta}
\frac{\Gamma\left(\frac{n+4-\beta}{2}\right)\Gamma\left(\frac{n+2}{2}\right)}{(n-\beta)\Gamma\left(\frac{\beta}{2}\right)}\|\nu\|^{\beta-n}
&\text{if $\beta\ne n$},\\
-\frac{2n}{\delta^2}\left(
2\log\|\nu\|+
\log\left(\frac{\delta^2}{4}\right)+\gamma-\psi(\frac{n}{2})\right)
&\text{if $\beta = n$},
\end{cases}
\end{equation*}
where $\gamma$ is Euler's constant and $\psi$ is the digamma function.
\end{theorem}
%\section{Temp title}
\subsection{Explicit formula for the extended  operator for $\beta <n+4$.}\label{sec:L-explicit}
The extended operator $\Ldel$ as described in   Section~\ref{sec:prelim}  is only defined indirectly through the Fourier transform as a pseudo-operator for the case $\beta>n+2$. In this section, we provide an explicit differential-integral formula for $\Ldel$ when $\beta<n+4$. We emphasize that such representation can be generalized in a similar way to $\beta>n+4$ with $\beta\ne n+6,n+8,\cdots$.
\begin{theorem}\label{thm:Laplace_extended}
Let $\beta <n+4$. Then the nonlocal Laplace operator defined as the pseudo-differential operator through its multipliers given in \eqref{eq:multiplier-general} can be represented explicitly in the form
\begin{equation}\label{eq:nonlocal_laplacian_extended}
    L^{\delta,\beta}u(x)=\Delta u(x)+c^{\delta,\beta}\int_{B_\delta(0)}\frac{u(x+z)-u(x)-\nabla\nabla  u(x):\frac{z\otimes z}{2}}{\|z\|^\beta}\;dz,
\end{equation}
for sufficiently differentiable fields $u$.
Moreover, its multipliers can be expressed as
\begin{equation}
    \label{eq:multiplier-cosine0-extended}
    \mdel(\nu)=-\|\nu\|^2+\cdel\int_{B_\delta(0)}\frac{\cos(\nu\cdot z)-1+\frac{(\nu\cdot z)^2}{2}}{\|z\|^\beta}\:dz.
\end{equation}
\end{theorem}
\begin{proof}
We begin with representing $u(x)$ through its Fourier transform as
\begin{align*}
    u(x)=\frac{1}{(2\pi)^n}\int_{\mathbb{R}^n}\hat{u}(\nu)e^{i \nu\cdot x}d\nu.
\end{align*}
Applying the operator $\Ldel$  gives
\begin{align*}
    L^{\delta,\beta}u(x)&=\Delta u(x)+c^{\delta,\beta}\int_{B_\delta(0)}\frac{u(x+z)-u(x)-\nabla\nabla  u(x):\frac{z\otimes z}{2}}{\|z\|^\beta}\;dz\\
    &=\frac{1}{(2\pi)^n}\int_{\mathbb{R}^n}\hat{u}(x)\Delta e^{i \nu\cdot x}d\nu\\ 
    &\qquad+\frac{\cdel}{(2\pi)^n}\int_{B_\delta(0)}\|z\|^{-\beta}\left(\int_{\mathbb{R}^n}\hat{u}(\nu) e^{i \nu\cdot (x+z)}d\nu-\int_{\mathbb{R}^n}\hat{u}(\nu) e^{i \nu\cdot x}d\nu\right)dz\\
    &\qquad-\frac{\cdel}{(2\pi)^n}\int_{B_\delta(0)}\|z\|^{-\beta}\left(\int_{\mathbb{R}^n}\hat{u}(\nu) \nabla\nabla e^{i \nu\cdot x}d\nu\right):\frac{z\otimes z}{2}dz\\
    &=\frac{1}{(2\pi)^n}\int_{\mathbb{R}^n}\hat{u}(x)(-\|\nu\|^2) e^{i \nu\cdot x}d\nu \\
    &\qquad+\frac{\cdel}{(2\pi)^n}\int_{B_\delta(0)}\|z\|^{-\beta}\int_{\mathbb{R}^n}\hat{u}(\nu) e^{i \nu\cdot x}\left(e^{i \nu\cdot z}-1\right) d\nu dz\\
    &\qquad-\frac{\cdel}{(2\pi)^n}\int_{B_\delta(0)}\|z\|^{-\beta}\left(\int_{\mathbb{R}^n}\hat{u}(\nu)\left(-\nu \otimes \nu\right) e^{i \nu\cdot x}d\nu\right):\frac{z\otimes z}{2}dz\\
    &=\frac{1}{(2\pi)^n}\int_{\mathbb{R}^n}\left[-\|\nu\|^2+\cdel\int_{B_\delta(0)}\frac{\cos(\nu\cdot z)-1+\frac{(\nu\cdot z)^2}{2}}{\|z\|^\beta}\:dz\right]\hat{u}(\nu)e^{i \nu \cdot x}d\nu\\
    &:=\frac{1}{(2\pi)^n}\int_{\mathbb{R}^n}\mdel(\nu)\hat{u}(\nu)e^{i \nu \cdot x}d\nu.
\end{align*}
Moreover, using the spherical coordinates $\nu\cdot z=\|\nu\|r\cos{\phi_1}$, we obtain
\begin{align*}
 \!\!\!\!   \mdel(\nu)&=-\|\nu\|^2+\cdel\int_{B_\delta(0)}\frac{\cos(\nu\cdot z)-1+\frac{(\nu\cdot z)^2}{2}}{\|z\|^\beta}\:dz\\
    &=-\|\nu\|^2+\cdel\int_{0}^{\delta}r^{n-1-\beta}\int_{S^{n-1}}\left[\cos(\|\nu\|r\cos{\phi_1})-1\right.\\
    &\qquad\qquad\qquad\qquad\left.+\|\nu\|^2r^2\cos^2\phi_1/2\right]d_{S^{n-1}}Vdr\\
    &=-\|\nu\|^2+\cdel\int_{0}^{\delta}r^{n-1-\beta}\int_{S^{n-1}}\left[\cos(\|\nu\|r\cos{\phi_1})-1\right]d_{S^{n-1}}Vdr\\
    &\qquad+\frac{\cdel\|\nu\|^2}{2}\int_{0}^{\delta}r^{n+1+\beta}\int_{0}^{\pi}\cos^2\phi_1\sin^{n-2}\phi_1d\phi_1S_{n-2}dr.
\end{align*}
Using ~\cite[equation (15)]{alali2021fourier} we have
\begin{align*}
    \mdel(\nu)&=-\|\nu\|^2+\cdel\int_{0}^{\delta}r^{n-1-\beta}2\pi^{n/2}\sum_{k=1}^{\infty}\frac{\left(-\frac{1}{4}\|\nu\|^2r^2\right)^k}{k!\Gamma(\frac{n}{2}+k)} dr\\
    &\qquad+\frac{\cdel\|\nu\|^2}{2}\int_{0}^{\delta}r^{n+1-\beta}\sqrt{\pi}\left[\frac{\Gamma(\frac{n-1}{2})}{\Gamma(\frac{n}{2})}-\frac{\Gamma(\frac{n+1}{2})}{\Gamma(\frac{n+2}{2})}\right]\frac{2\pi^{(n-1)/2}}{\Gamma(\frac{n-1}{2})}dr\\
    &=-\|\nu\|^2+2\pi^{n/2}\cdel\sum_{k=1}^{\infty}\frac{\left(-\frac{1}{4}\|\nu\|^2\right)^k}{k!\Gamma(\frac{n}{2}+k)}\frac{\delta^{n-\beta+2k}}{n-\beta+2k}
    +\cdel\|\nu\|^2\frac{\delta^{n+2-\beta}}{n+2-\beta}\frac{(n-1)\pi^{n/2}}{\Gamma(\frac{n+2}{2})}\\
    &=2\pi^{n/2}\cdel\sum_{k=1}^{\infty}\frac{\left(-\frac{1}{4}\|\nu\|^2\right)^k}{k!\Gamma(\frac{n}{2}+k)}\frac{\delta^{n-\beta+2k}}{n-\beta+2k}.
\end{align*}
The series in the last equality was shown in \cite{alali2021fourier} to coincide with the hypergeometric function given in \eqref{eq:multiplier-general}, thus showing the result.

\end{proof}

\begin{proposition}
Let $\beta <n+2$. 
%and assume that $u$ is sufficiently differentiable.
Then, the definition of $\Ldel u$ given in \eqref{eq:nonlocal_laplacian_extended} coincides with that in \eqref{eq:nonlocal_laplacian}. Moreover, the multipliers given in \eqref{eq:multiplier-cosine0-extended} coincide with that in \eqref{eq:multiplier-cosine0}.
\end{proposition}
\begin{proof}
This is a direct application of the following observation. For $\beta<n+2$ and sufficiently differentiable $u$, we have 
\begin{align*}
    \Delta u(x)=\frac{\cdel}{2}\int_{B_\delta(0)}\frac{z \otimes z}{\|z\|^\beta}\;dz:\nabla \nabla u(x).
\end{align*}
On the other hand, for $\nu\in\mathbb{R}^n$
\begin{align*}
    \|\nu\|^2=I:(\nu\otimes\nu)=\left[\frac{\cdel}{2}\int_{B_\delta(0)}\frac{z\otimes z}{\|z\|^\beta}dz\right]:(\nu\otimes\nu)=\frac{\cdel}{2}\int_{B_{\delta}(0)}\frac{(\nu\cdot z)^2}{\|z\|^\beta}dz.
\end{align*}

\end{proof}

In \cite{alali2021fourier}, the convergence of the original operator $\Ldel$ for $\beta<n+2$ in \eqref{eq:multiplier-general}to the Laplacian operator $\Delta$ under two limits as $\delta\to 0^+$ or as$\beta\to (n+2)^-$ has been established. We now generalize this result to the  extended operator $\Ldel$ in \eqref{eq:nonlocal_laplacian_extended} under the assumption that $\beta<n+3$. A similar result can be obtained when $n+3\le \beta<n+4$ in requirement of more regularity imposed on $u(x)$.

\begin{theorem}
    For $u\in C^3(\mathbb{R}^n)$, we have
    \begin{align*}
        \lim_{\delta\to 0^+}\Ldel u(x)=\Delta u(x)\quad \text{for }\beta<n+3,
    \end{align*}
    and
    \begin{align*}
        \lim_{\beta\to n+2}\Ldel u(x)=\Delta u(x)\quad \text{for }\delta>0.
    \end{align*}
\end{theorem}

\begin{proof}
    Verifying that \eqref{eq:nonlocal_laplacian_extended} is well-defined for $\beta<n+3$ and $u\in C^3(\mathbb{R}^n)$ is straightforward. Thanks to the Taylor's theorem, we get the following expansion of $u$ around $x$:
    \begin{align*}
        u(x+z)=u(x)+\nabla u(x):z+\frac{1}{2}\nabla\nabla u(x):z\otimes z+ R(u;x,z),
    \end{align*}
    where the remainder $R(u;x,z)$ is of the form
    \begin{align*}
        R(u;x,z)=\frac{1}{6}\nabla\nabla\nabla u(x+sz):z\otimes z\otimes z,
    \end{align*}
    for some scalar $s\in [0,1]$. Observe that \eqref{eq:nonlocal_laplacian_extended} now becomes
    \begin{align*}
        \Ldel u(x)-\Delta u(x)=\cdel\nabla u(x):\int_{B_\delta(0)}\frac{z}{\|z\|^\beta}dz+\cdel \int_{B_\delta(0)}\frac{R(u;x,z)}{\|z\|^\beta}dz.
    \end{align*}
    The first expression on the right hand side is identically zero due to symmetry. Hence, this shows that 
    \begin{align*}
        |\Ldel u(x)-\Delta u(x)|=\left|\cdel \int_{B_\delta(0)}\frac{R(u;x,z)}{\|z\|^\beta}dz\right|\le H_x|\cdel| \int_{B_\delta(0)}\frac{\|z\|^3}{\|z\|^\beta}dz.
    \end{align*}
    Here $H_x$ denotes 
    \begin{align*}
        H_x:=\frac{1}{6}\max\limits_{\substack{s\in [0,1]\\z\in \overline{B_\delta(0)}}}\left|\nabla\nabla \nabla u(x+sz)\right|.
    \end{align*}
    A simple calculation yields that 
    \begin{align*}
        \int_{B_\delta(0)}\frac{\|z\|^3}{\|z\|^\beta}dz=\frac{2\pi^{n/2}}{\Gamma\left(\frac{n}{2}\right)}\frac{\delta^{3+n-\beta}}{3+n-\beta}
    \end{align*}
    This observation indicates that
    \begin{align*}
        |\Ldel u(x)-\Delta u(x)|\le H_x\frac{2n|2+n-\beta|}{3+n-\beta}\delta,
    \end{align*}
    from which the convergence  follows.
\end{proof}

\subsection{Monotonicity of the multipliers}\label{sec:monotonicity}
In this section, we include a result on the monotonicity of the multipliers as a function of $\beta$. 
This lemma will be essential to prove convergence of solutions of the nonlocal wave equation to the solution of the classical wave equation in the limit as $\beta\rightarrow n+2$, as given in Theorem~\eqref{thm:converge1}  and Theorem~\eqref{converge2}.
%For the next theorem, we will make use of the following lemma on the monotonicity of the multipliers in the kernel parameter $\beta$. 
%which was proved in \cite{alali2021fourier}.

\begin{lemma}
\label{lem:monot1}
Fixing $\delta>0$ and $\nu \neq 0$, the function $\beta\mapsto m^{\delta,\beta}(\nu)$ is monotonically decreasing for $\beta \in \left(\frac{n+4}{2}, n+4\right)$.
\end{lemma}
\noindent

\begin{proof}
%[Proof of Lemma~\ref{lem:monot1}]
We express the multipliers as
\begin{equation*}
	m^{\delta,\beta}(\nu) = -\|\nu\|^2\,_2F_3\left(1,a;2,a+1,b;-\frac{1}{4}\|\nu\|^2\delta^2\right),
\end{equation*}
 where
\begin{equation*}
a = \frac{n+2-\beta}{2}\quad\text{and}\quad b=\frac{n+2}{2}.
\end{equation*}
For $z< 0$ define the function
\begin{equation*}
\phi(a) := {}_2F_3(1,a;2,a+1,b;z).
\end{equation*}
 Hence, if we let $z=-\frac{1}{4}\delta^2\|\nu\|^2$, then
\begin{equation*}
m^{\delta,\beta}(\nu) = -\|\nu\|^2\phi(a).
\end{equation*}
Note that if $\phi(a)$ is strictly monotonically decreasing in $a$ and $\nu\ne 0$, then $m^{\delta,\beta}(\nu)$ is strictly monotonically decreasing in $\beta$.

We begin with an application of~\cite{NIST:DLMF} \href{https://dlmf.nist.gov/16.5#E2}{(16.5.2)}, which provides the representation
\begin{equation*}
_2F_3\left(1,a;2,a+1,b;z\right) = 
\int_0^1{}_1F_2(a;a+1,b;zt)\;dt.
\end{equation*}
\noindent
Thus, a sufficient condition for the monotonicity of $\phi$ is that the function
\begin{equation*}
\psi(a) := {}_1F_2(a;a+1,b;z)
\end{equation*}
is strictly monotonically decreasing for any $z<0$.  In general, differentiating hypergeometric functions with respect to their parameters is messy, but the special structure of the arguments in this case makes it relatively straightforward. Written as a series,
\begin{equation*}
\psi(a) = \sum_{k=0}^\infty \frac{(a)_k}{(a+1)_k(b)_k}\frac{z^k}{k!}
= \sum_{k=0}^\infty \frac{a}{a+k} \frac{1}{(b)_k}\frac{z^k}{k!}.
\end{equation*}
According to the Weierstrass M-test, the series converges uniformly. Thus,
\begin{equation*}
\begin{split}
\psi'(a) &= \sum_{k=1}^\infty \frac{k}{(a+k)^2}\frac{1}{(b)_k}\frac{z^k}{k!}
= \sum_{k=0}^\infty \frac{1}{(a+1+k)^2}\frac{1}{(b)_{k+1}}\frac{z^{k+1}}{k!}\\
&= \frac{z}{b(a+1)^2}\sum_{k=0}^\infty\frac{(a+1)^2}{(a+1+k)^2}\frac{1}{(b+1)_k}\frac{z^k}{k!}\\
&= \frac{z}{b(a+1)^2}\sum_{k=0}^\infty\frac{(a+1)_k(a+1)_k}{(a+2)_k(a+2)_k(b+1)_k}\frac{z^k}{k!}\\
&= \frac{z}{b(a+1)^2}\;{}_2F_3(a+1,a+1;a+2,a+2,b+1;z).
\end{split}
\end{equation*}
Since $z<0$, if we show that $_2F_3(a+1,a+1;a+2,a+2,b+1;z)>0$, it follows that $\psi'(a)<0$.
\noindent
Now consider~\cite{cho2021newton}.  As described in the introduction, of that paper, we may start from the identity that
\begin{equation*}
\begin{split}
_2F_3\left(a+1,a+1;a+1,a+\frac{3}{2},2(a+1);-x^2\right) &=
{}_1F_2\left(a+1;a+\frac{3}{2},2(a+1);-x^2\right) \\
&= \frac{\Gamma\left(a+\frac{3}{2}\right)^2}{(x/2)^{2a+1}}
J_{a+\frac{1}{2}}(x)^2 \ge 0.
\end{split}
\end{equation*}
Following that, we may use the transference principle (see~\cite[Prop.~2.1]{cho2021newton}) to assert that, since $a+2>a+3/2>a+1$ and $b+1>2(a+1)$ (as long as $\beta > \frac{n+4}{2}$), it follows that $\psi'(a)<0$ as desired.
\end{proof}

\begin{remark}
We note that a version of Lemma~\eqref{lem:monot1} has appeared in \cite{alali2021fourier}, however the proof contained an error. This is fixed in Lemma~\eqref{lem:monot1} above, whose proof is based on the hypergeometric representation \eqref{eq:multiplier-general}. In addition, a simpler proof based on the representation given by \eqref{eq:multiplier-cosine0-extended}, but for a more restrictive condition on $\beta$, is included below for clarity of the presentation.
\end{remark}
\begin{lemma}
\label{lem:beta}
Let $n+2\le \beta<n+4$ and assume $\beta'< \beta$. 
%$\beta'<\beta<n+4$. 
Then, for all $\nu\ne 0$,
\[
\mdel(\nu)<m^{\delta,\beta'}(\nu)<0.
\]
\end{lemma}
\begin{proof}
Obviously, when $\beta<n+2$, $\cdel>0$. From \eqref{eq:multiplier-cosine0} we have $\mdel(\nu)<0$. Otherwise, when $n+2\le \beta< n+4$, $\cdel \le 0$. From \eqref{eq:multiplier-cosine0-extended}, we also get $\mdel(\nu)<0$.\\
\begin{itemize}
    \item Case $\beta'<n+2\le \beta$: We have $c^{\delta,\beta'}>0\ge \cdel$. This leads to
\begin{align*}
    \mdel(\nu)\le -\|\nu\|^2< m^{\delta,\beta'}(\nu).
\end{align*}

    \item Case $n+2\le \beta'<\beta <n+4$: By changing the variable $z=\delta w$
\begin{align*}
    \mdel(\nu)&=-\|\nu\|^2+\cdel\delta^{n-\beta}\int_{B_1(0)}\frac{\cos(\nu\cdot \delta w)-1+\frac{(\nu\cdot \delta w)^2}{2}}{\|w\|^\beta}\:dw\\
    &=-\|\nu\|^2-\frac{2(\beta-(n+2))\Gamma\left(\frac{n}{2}+1\right)}
    {\pi^{n/2}\delta^{2}}\int_{B_1(0)}\frac{\cos(\nu\cdot \delta w)-1+\frac{(\nu\cdot \delta w)^2}{2}}{\|w\|^\beta}\:dw.
\end{align*}
Notice that when $n+2\le\beta'<\beta<n+4$, we have $\beta-(n+2)>\beta'-(n+2)\ge 0$. Also, $\|w\|^\beta<\|w\|^{\beta'}$ and $\cos(\nu\cdot \delta w)-1+\frac{(\nu\cdot \delta w)^2}{2}\ge 0$ for $w\in B_1(0)$. All of these yield
\begin{align*}
    \mdel(\nu)< m^{\delta,\beta'}(\nu).
\end{align*}

%    \item Case $\beta'<\beta<n+2:$
\end{itemize}

\end{proof}

\subsection{The periodic torus and the eigenvalues of $\Ldel$.}
Define the periodic torus as
\begin{equation*}
\T^n = \prod_{j=1}^n[0,2\pi],\qquad \text{with}\quad j=1,2,\ldots,n.
\end{equation*}
\noindent
The functions $\{e^{i k\cdot x}\}_{k\in\mathbb{Z}^n}$ form a complete set in $L^2(\T^n)$.  Moreover, when $\beta<n+2$,
\begin{equation*}\label{eq:phi}
  \Ldel e^{i k\cdot x} =
  \left(
  \cdel\int_{B_\delta(0)}\frac{e^{i k\cdot z}-1}{\|z\|^\beta}\;dz
  \right)e^{i k\cdot x}
  = m^{\delta,\beta}(k)e^{i k\cdot x},
\end{equation*}
which implies that $e^{i k\cdot x}$ is an eigenfunction of $\Ldel$ with eigenvalue $m^{\delta,\beta}(k)$.

In fact, the above observation on the eigenvalues of $\Ldel$ for $\beta<n+2$,  carries over to the  extended operator with any $\beta\ge n+2$ with $\beta\ne n+4,n+6,\cdots$. To see this, we give a proof  for the case when $\beta<n+4$.
\begin{proposition}
Let $\beta<n+4$. Then, for $k \in \mathbb{Z}^n$, $(\mdel(k),e^{ik\cdot x})$ is an eigenpair for the operator defined in \eqref{eq:nonlocal_laplacian_extended}.
\end{proposition}
\begin{proof}
We observe that
\begin{align*}
    \Ldel e^{i k\cdot x}&=\Delta e^{i k\cdot x}+\cdel\int_{B_\delta(0)}\frac{e^{i k\cdot (x+z)}-e^{i k\cdot x}-\nabla \nabla e^{i k\cdot x}:\frac{z\otimes z}{2}}{\|z\|^\beta}\:dz\\
    &=-\|k\|^2e^{i k\cdot x}+\cdel\int_{B_\delta(0)}\frac{(e^{i k\cdot z}-1)e^{i k\cdot x}+k\otimes k e^{i k\cdot x}:\frac{z\otimes z}{2}}{\|z\|^\beta}\;dz\\
    &=\left(-\|k\|^2+\cdel\int_{B_\delta(0)}\frac{\cos(k\cdot z)-1+\frac{(k\cdot z)^2}{2}}{\|z\|^\beta}\:dz\right)e^{i k \cdot x}=\mdel(k)e^{i k\cdot x}.
\end{align*}
\end{proof}

Throughout the remaining part of this article, and to simplify notation,
%To simplify the presentation, 
we use  $m_k$ to denote  $m^{\delta,\beta}(k)$ for any $k\in\mathbb{Z}^n$. We observe that for the case $\beta=n+2$, the eigenvalue $m_k$ is equal to $-\|k\|^2$.

\subsection{Sobolev spaces on the torus}

\noindent
  For $q\in\mathbb{R}$, let $H^q(\T^n)$ denote the space of all periodic distributions $g$ on $\mathbf{T}^n$ such that
\begin{equation*}
    \sum_{k\in\mathbb{Z}^n}(1+\|k\|^2)^q|\hat{g}_k|^2
    < \infty.
\end{equation*}
Define the norm on $H^q(\T^n)$ by
\begin{align*}
    \|g\|_{H^q(\T^n)}=\left(\sum_{k\in\mathbb{Z}^n}(1+\|k\|^2)^q|\hat{g}_k|^2\right)^{1/2}<\infty,
\end{align*}
for any $g\in H^q(\T^n)$.

\begin{lemma}
Let $f\in H^q(\T^n)$. For any nonzero $k\in \mathbb{Z}^n$, there exists $C>0$ such that
\begin{equation}\label{Forier-decay} 
    |\hat{f}_k|\le \frac{C}{\|k\|^q}.
\end{equation}
\end{lemma}
\begin{proof}
For any $k\ne 0$ in $\mathbb{Z}^n$, there is a constant $D>0$ satisfying
\begin{align*}
    (1+\|k\|^2)^q\ge D\|k\|^{2s}.
\end{align*}
Thus, using the fact that $f$ is in $H^q(\T^n)$, we obtain
\[
\infty>\sum_{0\ne k\in\mathbb{Z}^n}(1+\|k\|^2)^q |\hat{f}_k|^2\ge D \sum_{0\ne k\in\mathbb{Z}^n}\|k\|^{2q} |\hat{f}_k|^2.
\]
This implies that $\|k\|^{2q} |\hat{f}_k|^2\rightarrow 0$ as $\|k\|\rightarrow \infty$, from which the result follows.
%Since $f$ is in $H^s(\T^n)$, we have
%\begin{align*}
%    \infty>\|f\|^2_{H^s(\T^n)}\ge D|\hat{f}_k|(1+\|k\|^2)^{s/2}\ge |\hat{f}_k|\|k\|^s.
%\end{align*}
%This leads us to the result of the lemma.
\end{proof}

In  subsequent sections, we consider maps from the real line to periodic distributions and use the following notion of derivative for these maps.
\noindent
\begin{definition}
Suppose $T:\left[0,\infty\right)\to H^q(\T^n)$ maps $t\ge 0$ to a periodic distribution $T(t)$ in $H^q(\T^n)$. Define the derivative of $T$ at $t$ in the sense of Gateaux derivative
\begin{align*}
    T^\prime(t)=\frac{d}{dt}T(t)=\lim_{h\to 0}\frac{T(t+h)-T(t)}{h},
\end{align*}
where the limit is with respect to the norm in $H^q(\T^n)$.\\
\end{definition}
Let $p\in\mathbb{N}$ and $q\in\mathbb{R}$, and  denote by $C^p([0,\infty), H^q(\T^n))$ the space of all maps from $[0,\infty)$ to $H^q(\T^n)$ that is $p$ times differentiable. 
\begin{lemma}\label{lem:derivative}
Given a map $T:\left[0,\infty\right)\to H^q(\T^n)$ and suppose that the derivative $T^\prime(t)$ exists. Then for any $k\in\mathbb{Z}^n$, the Fourier coefficient $\hat{T}^\prime_k(t)$ of $T^\prime(t)$ is the derivative, in the classical sense, of the Fourier coefficient $\hat{T}_k(t)$. 
\end{lemma}
\begin{proof}
Represent $T(t)$ and $T^\prime(t)$ as
\begin{align*}
    T(t)=\sum_{k\in\mathbb{Z}^n}\hat{T}_k(t)e^{i k\cdot x},
\end{align*}
and
\begin{align*}
    T^\prime(t)=\sum_{k\in\mathbb{Z}^n}\hat{T}^\prime_k(t)e^{i k\cdot x}.
\end{align*}
From the derivative definition, we have
\begin{align*}
    0&=\lim_{h\to 0}\left\|\frac{T(t+h)-T(t)}{h}-T^\prime(t)\right\|_{H^q(\T^n)}\\
    &=\lim_{h\to 0}\sum_{k\in\mathbb{Z}^n}(1+\|k\|^2)^q\left|\frac{\hat{T}_k(t+h)-\hat{T}_k(t)}{h}-\hat{T}^\prime_k(t)\right|^2.
\end{align*}
This proves that
\begin{align*}
    \lim_{h\to 0}\frac{\hat{T}_k(t+h)-\hat{T}_k(t)}{h}=\hat{T}^\prime_k(t),
\end{align*}
for any $k\in \mathbb{Z}^n$, from which the result follows.
%. This is the definition for derivative in the classical sense.
\end{proof}

\begin{definition}
Suppose $T:\left[0,\infty\right)\to H^q(\T^n)$ is such that for any $t\ge 0$, $T(t)$  has the form
\begin{align*}
    T(t)=\sum_{k\in\mathbb{Z}^n}\hat{T}_k(t)e^{ik\cdot x}.
\end{align*}
Then, define the map $\Ldel T:\left[0,\infty\right)\to H^q(\T^n)$ through the representation
\begin{align*}
    \Ldel T(t)=\sum_{k\in\mathbb{Z}^n}m_k\hat{T}_k(t)e^{i k\cdot x},
\end{align*}
for any $t\ge 0$.
%The operation Next we provide a  definition for the  The following definition provides a meaning for the nonlocal Laplace operator $\Ldel$
\end{definition}
In particular, when $q<0$, the above notion extends the definition of the nonlocal Laplacian $\Ldel$ to periodic distributions. We emphasize that the discussion and results here apply for the extended nonlocal Laplace operator.

\section{Homogeneous nonlocal wave equation}\label{sec:Homogeneous}
In this section, we focus on the following nonlocal wave equation with initial data and zero forcing term
\begin{align}\label{eq:Homo-wave-general}
    \begin{cases}
    u_{tt}(x,t)&=\Ldel u(x,t),\quad x\in\T^n, \; t>0,\\
    u(x,0)&=f(x),\\
    u_t(x,0)&=g(x).
    \end{cases}
\end{align}
We emphasize  that the results of this section hold for any nonlocal operator $\Ldel$ with $\beta<n+4$. We recall from Sections \ref{sec:prelim} and \ref{sec:L-explicit} that $\beta<n$ corresponds  to an integrable kernel, $n\le \beta<n+2$ corresponds  to a singular kernel, $\beta=n+2$ corresponds  to $\Ldel=\Delta$ (the Laplacian), and $n+2<\beta<n+4$ corresponds  to an integro-differential operator as given in \eqref{eq:nonlocal_laplacian_extended}.  

In order to study the existence, uniqueness, and regularity of solutions to \eqref{eq:Homo-wave-general} over the space of periodic distributions, we consider the identification $U(t)=u(\cdot,t)$, with $U:[0,\infty)\rightarrow H^q(\T^n)$, for some $q\in\mathbb{R}$. In the following sections, spatial regularity refers to $U(t)$ being in $H^q(\T^n)$, for some $q\in\mathbb{R}$, for a fixed $t\ge 0$. Temporal regularity refers to the differentiability of the map $U$ as a function of time $t$.
\subsection{Spatial and temporal regularity over periodic distributions}
Fix $f\in H^{s_1}(\T^n)$ and $g\in H^{s_2}(\T^n)$, with $s_1$ and $s_2$ in $\mathbb{R}$. For any $t\ge 0$, define
\begin{align}\label{EqU-homo}
\nonumber
    U(t)&:=\sum_{k\in \mathbb{Z}^n}\hat{U}_k(t)e^{ik\cdot x}\\
    &=(\hat{f}_0+\hat{g}_0t)+\sum_{0\ne k\in \mathbb{Z}^n}\left[\hat{f}_k\cos(\sqrt{-m_k}t)+\frac{\hat{g}_k}{\sqrt{-m_k}}\sin(\sqrt{-m_k}t)\right]e^{i k\cdot x}.
\end{align}
We will show in Theorem~\eqref{thm:existence} that $U$ is the unique solution to the nonlocal wave equation. We begin by studying the spatial and temporal differentiability of $U$.
\begin{theorem}\label{thm:spatial_reg1}
Let $n\ge 1$, $\delta>0$, and $\beta<n+4$. For any $t\ge 0$, $U(t)\in H^s(\T^n)$, where $s=\min\{s_1,s_2+\theta\}$ with $\theta=\max\{0,\frac{\beta-n}{2}\}$.
\end{theorem}
\begin{proof}
Observe that
\begin{align*}
  \!\!\!\!\!\!  \sum_{k\in\mathbb{Z}^n}(1+\|k\|^2)^s|\hat{U}_k(t)|^2&\le \sum_{0\ne k\in \mathbb{Z}^n}(1+\|k\|^2)^s|\hat{f}_k|^2
    %&\qquad
    +\sum_{0\ne k\in \mathbb{Z}^n}(1+\|k\|^2)^s\frac{|\hat{g}_k|^2}{|m_k|}+|\hat{f}_0|^2+|\hat{g}_0|^2t^2.
\end{align*}
Since $f\in H^{s_1}$ and $s\le s_1$, we get
\begin{align*}
    \sum_{0\ne k\in \mathbb{Z}^n}(1+\|k\|^2)^s|\hat{f}_k|^2\le \sum_{0\ne k\in \mathbb{Z}^n}(1+\|k\|^2)^{s_1}|\hat{f}_k|^2<\infty.
\end{align*}
Additionally, notice that
\begin{align*}
    \sum_{0\ne k\in \mathbb{Z}^n}(1+\|k\|^2)^s\frac{|\hat{g}_k|^2}{|m_k|}=\sum_{0\ne k\in \mathbb{Z}^n}\frac{(1+\|k\|^2)^{s-s_2}}{|m_k|}(1+\|k\|^2)^{s_2}|\hat{g}_k|^2.
\end{align*}
Hence, the proof will be complete provided that
\begin{align*}
    \frac{(1+\|k\|^2)^{s-s_2}}{|m_k|},
\end{align*}
is bounded for large values of $\|k\|$.
To show this, we consider two cases. First, when $\beta\le n$, then $s\le s_2$ and by Theorem ~\eqref{thm:asymptotics} there exists $r_1>0$ and $C_1>0$ such that   $|m_k|\ge C_1$ for $\|k\|\ge r_1$. This leads  to
\begin{align*}
    \frac{(1+\|k\|^2)^{s-s_2}}{|m_k|}\le \frac{1}{C_1}.
\end{align*}
Similarly, for $\beta> n$, then $s \le s_2+\frac{\beta-n}{2}$ and by Theorem \eqref{thm:asymptotics} there exists $r_2>0$ and $C_2>0$ such that  $|m_k|\ge C_2\|k\|^{\beta-n}$, for $\|k\|\ge r_2$. This leads  to
\begin{align*}
    \frac{(1+\|k\|^2)^{s-s_2}}{|m_k|}\le \frac{1}{C_2}\left(\frac{1+\|k\|^2}{\|k\|^2}\right)^{\frac{\beta-n}{2}},
\end{align*}
which is bounded.
\end{proof}

\begin{theorem}\label{thm:regularity}
Let $n\ge 1$, $\delta>0$, and $\beta<n+4$. Let $U$ be the map given by \eqref{EqU-homo} with $s_1, s_2\in \mathbb{R}$, and $s=\min\{s_1,s_2+\theta\}$ with $\theta=\max\{0,\frac{\beta-n}{2}\}$. Thus,
\begin{enumerate}
    \item if $\beta<n$, then $U\in C^\infty([0,\infty),H^q(\T^n))$, for any $q\le s$,
    \item if $\beta=n$, then $U\in C^\infty([0,\infty),H^q(\T^n))$,  for any $q<s$, and
    \item if $\beta>n$, then $U\in C^{p+1}([0,\infty),H^q(\T^n))$, for any $q\in\mathbb{R}$ and any positive integer $p$    satisfying 
    \begin{align*}
        q-s_1+(p+1)\frac{\beta-n}{2}\le 0\quad\text{and}\quad q-s_2+p\frac{\beta-n}{2}\le 0.
    \end{align*}

\end{enumerate}
\end{theorem}

\begin{proof}
We will show this result using induction. Suppose $U$ is already differentiable up to $p$ times, then by Lemma \eqref{lem:derivative}, we  have
\begin{align*}
    U^{(p)}(t)=\sum_{k\in\mathbb{Z}^n}\hat{U}^{(p)}_k(t)e^{ik\cdot x}
    .
\end{align*}
where $\hat{U}^{(p)}_k(t)$ is given by
\begin{align*}
    \hat{U}^{(p)}_k(t)=\hat{f}_k(\sqrt{-m_k})^p\cos\left(\sqrt{-m_k}t+\frac{p\pi}{2}\right)+\hat{g}_k(\sqrt{-m_k})^{p-1}\sin\left(\sqrt{-m_k}t+\frac{p\pi}{2}\right),
\end{align*}
 for $k\ne 0$ and $\hat{U}^{(1)}_0(t)=\hat{g}_0$ and  $\hat{U}^{(p)}_0(t)=0$ for all $p\ge 2 $. We will show that 
\begin{align*}
    U^{(p+1)}(t)=\sum_{k\in\mathbb{Z}^n}\hat{U}^{(p+1)}_k(t)e^{ik\cdot x}
    .
\end{align*}
Observe that using the mean value theorem,
\begin{align*}
    T_h:&=\left\|\frac{U^{(p)}(t+h)-U^{(p)}(t)}{h}-\sum_{k\in\mathbb{Z}^n}\hat{U}^{(p+1)}_k(t)e^{ik\cdot x}\right\|^2_{H^q(\T^n)}\\
    &=\sum_{k\in\mathbb{Z}^n}(1+\|k\|^2)^q\left|\frac{\hat{U}^{(p)}_k(t+h)-\hat{U}^{(p)}_k(t)}{h}-\hat{U}^{p+1}_k(t)\right|^2\\
    &\le\sum_{0\ne k\in\mathbb{Z}^n}(1+\|k\|^2)^q|\hat{f}_k|^2(-m_k)^{p+1}\left|\cos\left(\sqrt{-m_k}(t+\xi_{p,k})+\frac{(p+1)\pi}{2}\right)\right.\\
    &\qquad\qquad\qquad\left.-\cos\left(\sqrt{-m_k}t+\frac{(p+1)\pi}{2}\right)\right|^2\\
    &\qquad+\sum_{0\ne k\in\mathbb{Z}^n}(1+\|k\|^2)^q|\hat{g}_k|^2(-m_k)^{p}\left|\sin\left(\sqrt{-m_k}(t+\eta_{p,k})+\frac{(p+1)\pi}{2}\right)\right.\\
    &\qquad\qquad\qquad\left.-\sin\left(\sqrt{-m_k}t+\frac{(p+1)\pi}{2}\right)\right|^2,
\end{align*}
for some $\xi_{p,k},\eta_{p,k}\in(0,h)$. This means that as $h\to 0$, each term of the summation will approach 0. The proof will be complete if we are able to pass the limit inside the summation. This is possible provided that the following terms are uniformly bounded
\begin{align*} 
    a_k:=(1+\|k\|^2)^{q-s_1}(-m_k)^{p+1}\quad\text{and}\quad b_k:=(1+\|k\|^2)^{q-s_2}(-m_k)^p.
\end{align*}
To establish the bound for these terms, we will consider three cases when $\beta<n$, $\beta=n$ and $\beta>n$.\\
When $\beta<n$, from Theorem \eqref{thm:asymptotics}, $-m_k$ is bounded, and this  provides the uniform boundedness of $a_k$ and $b_k$ provided that $q\le \min \{s_1,s_2\}$.\\
When $\beta=n$, fix an $\epsilon>0$. From Theorem \eqref{thm:asymptotics}, there exists $C,r>0$ such that $-m_k\le C\log\|k\|\le C(1+\|k\|^2)^{\epsilon/2}$ whenever $\|k\|\ge r$. Then for large $\|k\|$,
\begin{align*}
    a_k\le C(1+\|k\|^2)^{q-s_1+(p+1)\epsilon/2}\quad\text{and}\quad b_k\le C(1+\|k\|^2)^{q-s_2+p\epsilon/2}.
\end{align*}
This yields the uniform boundedness of $a_k$ and $b_k$ as long as $q-s_1+(p+1)\epsilon/2\le 0$ and $q-s_2+p\epsilon/2\le 0$. So given $q<\min\{s_1,s_2\}$, we can always pick a sufficiently  small $\epsilon$ that satisfies these two inequalities.\\
When $\beta>n$, from Theorem \eqref{thm:asymptotics}, there is $C,r>0$ such that $-m_k\le C\|k\|^{\beta-n}\le C(1+\|k\|^2)^{(\beta-n)/2}$ for $\|k\|\ge r$. Then for large $\|k\|$,
\begin{align*}
    a_k\le C(1+\|k\|^2)^{q-s_1+(p+1)\frac{\beta-n}{2}}\quad\text{and}\quad b_k\le C(1+\|k\|^2)^{q-s_2+p\frac{\beta-n}{2}}.
\end{align*}
This shows that the uniform boundedness of $a_k,b_k$ is established given the conditions
\begin{align*}
    q-s_1+(p+1)\frac{\beta-n}{2}\le 0\quad\text{and}\quad q-s_2+p\frac{\beta-n}{2}\le 0.
\end{align*}
\end{proof}

\noindent The following theorem is a direct application of Theorem \eqref{thm:spatial_reg1} and Theorem \eqref{thm:regularity}. The assumptions are to guarantee that the second derivative exists.
\begin{theorem}
\label{thm:existence}
Let $n\ge 1$, $\delta>0$, and $\beta<n+4$. Let $f\in H^{s_1}(\T^n)$ and $g\in H^{s_2}(\T^n)$, with  $s_1,s_2\in \mathbb{R}$. Let $q\in\mathbb{R}$ be such that
\begin{enumerate}
    \item $q\le \min\{s_1,s_2\}$, in the case when $\beta<n$,
    \item $q<\min\{s_1,s_2\}$, in the case when  $\beta=n$, and
    \item $q-s_1+(\beta-n)\le 0$ and $q-s_2+(\beta-n)/2\le 0$, in the case when $\beta>n$.
\end{enumerate}
Then,  $U(t)$ defined by \eqref{EqU-homo} is in $H^q(\T^n)$, for $t\ge 0$. Moreover,  $U$ is the unique solution to the nonlocal wave equation
\begin{align}
\label{eq:Homo-wave}
    \begin{cases}
    \frac{d^2}{dt^2}U(t)&=\Ldel U(t),\;\; t>0\\
    U(0)&=f,\\
    \frac{d}{dt}U(0)&=g.
    \end{cases}
\end{align}
Furthermore, when $\beta\le n$, $U\in C^\infty([0,\infty); H^s(\T^n))$ where $s=\min\{s_1,s_2\}$,  and when \\$\beta>n$, 
$U\in C^2([0,\infty); H^s(\T^n))$ where $s=\min\{s_1,s_2+\frac{\beta-n}{2}\}$.

\end{theorem}
%\begin{theorem}
%\label{thm:existence}
%Let $f\in H^{s_1}(\T^n)$ and $g\in H^{s_2}(\T^n)$, with  $s_1,s_2\in \mathbb{R}$. 
%Then,  $U$ defined by \eqref{EqU}  is the unique solution the nonlocal wave equation
%\begin{align}\label{eq:Homo-wave}
%    \begin{cases}
%    \frac{d^2}{dt^2}U(t)&=\Ldel U(t),\;\; t>0\\
%    U(0)&=f,\\
%    \frac{d}{dt}U(0)&=g.
%    \end{cases}
%\end{align}
%Moreover, when $\beta\le n$, $U(t)\in C^\infty([0,\infty); H^s(\T^n))$ where $s=\min\{s_1,s_2\}$,  and when $\beta>n$, 
%$U(t)\in C^2([0,\infty); H^s(\T^n))$ where $s=\min\{s_1,s_2+\frac{\beta-n}{2}\}$.
%
%\end{theorem}

\subsection{Convergence to the classical wave equation solution}
\label{sec:convergence1}
In this section, we will provide some results on the convergence of the solutions of the nonlocal wave equation \eqref{eq:Homo-wave} to the  solution of the corresponding local wave equation for two types of  limits: as $\delta \to 0^+$, with $\beta< n+4$ is being fixed,  or as $\beta\to n+2$, with $\delta>0$ is being fixed. The  local wave equation is defined by
\begin{align}\label{eq:Homo-wave-local}
    \begin{cases}
    \frac{d^2}{dt^2}U^0(t)&=\Delta U^0(t),\;\; t>0\\
    U^0(0)&=f,\\
    \frac{d}{dt}U^0(0)&=g,
    \end{cases}
\end{align}
and its solution is given by 
\begin{align}\label{EqU-homo-local}
    U^0(t)&:=(\hat{f}_0+\hat{g}_0t)+\sum_{0\ne k\in \mathbb{Z}^n}\left[\hat{f}_k\cos(\|k\| t)+\frac{\hat{g}_k}{\|k\|}\sin(\|k\| t)\right]e^{i k\cdot x}.
\end{align}
To emphasize the dependence of the solution of the nonlocal wave equation \eqref{eq:Homo-wave} on the parameters $\delta$ and $\beta$, in the remaining part of this section we denote the solution $U$ by $U^{\delta,\beta}$.

\begin{theorem}
Let $n\ge 1$, $\delta>0$, and $\beta<n+4$. Suppose $f\in H^{s_1}(\T^n)$ and $g\in H^{s_2}(\T^n)$ for some $s_1,s_2\in \mathbb{R}$. Let $U^0(t)\in H^s(\T^n)$ and $U^{\delta,\beta}(t)\in H^{s'}(\T^n)$ be the solutions to the local and nonlocal wave equations \eqref{eq:Homo-wave-local} and \eqref{eq:Homo-wave}, respectively, where $s=\min\{s_1,s_2+1\}$ and $s'=\min\{s_1,s_2+\theta\}$ with $\theta=\max\{0,\frac{\beta-n}{2}\}$. Let $t\ge 0$. Then, when $\beta \le n+2$,
\begin{equation*}
    \lim_{\delta\to 0^+}U^{\delta,\beta}(t)=U^0(t) \quad \text{in} \quad H^{s'}(\T^n),
\end{equation*}
and when $n+2<\beta<n+4$
\begin{equation*}
    \lim_{\delta\to 0^+}U^{\delta,\beta}(t)=U^0(t) \quad \text{in} \quad H^{s}(\T^n).
\end{equation*}
\end{theorem}

\begin{proof} 
It is easy to check that $H^s(\T^n)\subset H^{s'}(\T^n)$ when $\beta \le n+2$ and $H^{s'}(\T^n)\subset H^{s}(\T^n)$ when $n+2<\beta< n+4$. Hence all the limits make sense. From \eqref{EqU-homo} and \eqref{EqU-homo-local}, the Fourier coefficients $\hat{U}_k^{\delta,\beta}(t)$ in the nonlocal case and $\hat{U}_k(t)$ in the classical case ($\beta=n+2$) is determined by
\begin{align*}
    \hat{U}_k^{\delta,\beta}(t)&=\hat{f}_k\cos(\sqrt{-m_k}t)+\frac{\hat{g}_k}{\sqrt{-m_k}}\sin(\sqrt{-m_k}t),\\
    \hat{U}^0_k(t)&=\hat{f}_k\cos(\|k\|t)+\frac{\hat{g}_k}{\|k\|}\sin(\|k\|t),
\end{align*}
for nonzero $k$ and $\hat{U}_0^{\delta,\beta}(t)=\hat{U}^0_0(t)=\hat{f}_0+\hat{g}_0 t$. From these observation, for $s_0=s'$ when $\beta\le n+2$ and $s_0=s$ when $n+2<\beta<n+4$, one can get
\begin{align*}
    \|U^{\delta,\beta}(t)-U^0(t) \|^2_{H^{s_0}(\T^n)} &= \sum_{0\ne k\in\mathbb{Z}^n}(1+\|k\|^2)^{s_0}\|\hat{U}^{\delta,\beta}_k(t) - \hat{U}^0_k(t)\|^2 \\
&\le \sum_{0\ne k\in\mathbb{Z}^n}(1+\|k\|^2)^{s_0}|\hat{f}_k|^2|\cos(\sqrt{-m_k}t)-\cos(\|k\|t)|^2\\
&+\sum_{0\ne k\in\mathbb{Z}^n}(1+\|k\|^2)^{s_0}|\hat{g}_k|^2\left|\frac{\sin(\sqrt{-m_k}t)}{\sqrt{-m_k}}-\frac{\sin(\|k\|t)}{\|k\|}\right|^2.
\end{align*}
Since $f\in H^{s_1}(\T^n)$ and $g\in H^{s_2}(\T^n)$, one can pass the limit as $\delta\to 0^+$ term by term inside the summation provided the uniform boundedness of
\begin{align*}
    (1+\|k\|^2)^{s_0-s_1}|\cos(\sqrt{-m_k}t)-\cos(\|k\|t)|^2
\end{align*}
and 
\[
    (1+\|k\|^2)^{s_0-s_2}\left|\frac{\sin(\sqrt{-m_k}t)}{\sqrt{-m_k}}-\frac{\sin(\|k\|t)}{\|k\|}\right|^2,
\]
for large $\|k\|$. It is obvious that given $s_0\le s_1$,
\begin{align*}
    (1+\|k\|^2)^{s_0-s_1}|\cos(\sqrt{-m_k}t)-\cos(\|k\|t)|^2\le 4.
\end{align*}
For the second term, one gets
\begin{align*}
    (1+\|k\|^2)^{s_0-s_2}\left|\frac{\sin(\sqrt{-m_k}t)}{\sqrt{-m_k}}-\frac{\sin(\|k\|t)}{\|k\|}\right|^2
    \le (1+\|k\|^2)^{s_0-s_2}\left(\frac{1}{\sqrt{-m_k}}+\frac{1}{\|k\|}\right)^2.
\end{align*}
We consider two cases. When $\beta\le n$, then $s_0-s_2\le 0$. From Theorem~\eqref{thm:asymptotics}, there exists $r_1>0$ and $C_1>0$ such that $-m_k\ge C_1$ for $\|k\|\ge r_1$. When $\beta>n$, then $s_0-s_2\le (\beta-n)/2$ and $s_0-s_2\le 1$. From Theorem~\eqref{thm:asymptotics}, there exists $r_2>0$ and $C_2>0$ such that $-m_k\ge C_2\|k\|^{\beta-n}$ for $\|k\|\ge r_2$. These results lead  to the uniform boundedness of the right hand side. 
Now, taking the limit as $\delta\to0^+$ inside the summation term-wise and applying the pointwise convergence of ~\eqref{col:delta} completes the proof.
\end{proof}

\begin{theorem}\label{thm:converge1}
Let $n\ge 1$, $\delta>0$, and $\beta<n+4$. Suppose $f\in H^{s_1}(\T^n)$ and $g\in H^{s_2}(\T^n)$ for some $s_1,s_2\in\mathbb{R}$. Let $U^0(t)\in H^s(\T^n)$ and $U^{\delta,\beta}(t)\in H^{s'}(\T^n)$ be the solutions to the local and nonlocal wave equations \eqref{eq:Homo-wave-local} and \eqref{eq:Homo-wave},  respectively, where $s=\min\{s_1,s_2+1\}$ and $s'=\min\{s_1,s_2+\theta\}$ with $\theta=\max\{0,\frac{\beta-n}{2}\}$. Then, for any $t\ge 0$ and for any $\varepsilon\in(0,2)$ and $\varepsilon<n/2$,
\begin{equation*}
    \lim_{\beta\to n+2}U^{\delta,\beta}(t)=U(t) \quad \text{in} \quad H^{s_0}(\T^n),
\end{equation*}
where $s_0=\min\left\{s_1,s_2+\frac{2-\varepsilon}{2}\right\}$.
\end{theorem}
\begin{proof}
Fix an $\varepsilon \in (0,2),\varepsilon<n/2$ and define $\beta'=n+2-\varepsilon>\frac{n+4}{2}$. 
For simplicity, we would still denote 
\begin{align*}
    m_k:=m^{\delta,\beta}(k)\quad\text{and}\quad m'_k:=m^{\delta,\beta'}(k).
\end{align*}
Then for any $\beta $ such that $|\beta-(n+2)|<\varepsilon$ or $\frac{n+4}{2}<\beta'<\beta<n+2+\varepsilon<n+4$, we have $\frac{\beta-n}{2}>\frac{2-\varepsilon}{2}$. This shows that $s, s'>s_0$, so the limit makes sense.

\begin{align*}
    \|U^{\delta,\beta}(t)-U^0(t) \|^2_{H^{s_0}(\T^n)} &\le \sum_{0\ne k\in\mathbb{Z}^n}(1+\|k\|^2)^{s_0}\|\hat{f}_k\|^2|\cos(\sqrt{-m_k}t)-\cos(\|k\|t)|^2\\
&+\sum_{0\ne k\in\mathbb{Z}^n}(1+\|k\|^2)^{s_0}\|\hat{g}_k\|^2\left|\frac{\sin(\sqrt{-m_k}t)}{\sqrt{-m_k}}-\frac{\sin(\|k\|t)}{\|k\|}\right|^2.
\end{align*}
Similar to the previous theorem, we need to establish the uniform boundedness of
\begin{align*}
    (1+\|k\|^2)^{s_0-s_1}|\cos(\sqrt{-m_k}t)-\cos(\|k\|t)|^2\\
\end{align*}
and
\[
    (1+\|k\|^2)^{s_0-s_2}\left|\frac{\sin(\sqrt{-m_k}t)}{\sqrt{-m_k}}-\frac{\sin(\|k\|t)}{\|k\|}\right|^2,
\]
for large $\|k\|$. Again, the first term is always bounded by 4 given $s_0\le s_1$. For the second term, from Theorem~\eqref{thm:asymptotics} there exists $r>0$ and $C,D>0$ such that $\sqrt{-m'_k}\ge C(1+\|k\|^2)^{\frac{2-\varepsilon}{4}}$ and $\|k\|\ge D(1+\|k\|^2)^{1/2}$ for $\|k\|\ge r$. This implies that
\begin{align*}
    &(1+\|k\|^2)^{s_0-s_2}\left|\frac{\sin(\sqrt{-m_k}t)}{\sqrt{-m_k}}-\frac{\sin(\|k\|t)}{\|k\|}\right|^2\\ 
    &\le (1+\|k\|^2)^{\frac{2-\varepsilon}{2}}\left(\frac{1}{\sqrt{-m_k}}+\frac{1}{\|k\|}\right)^2\\
    &\le (1+\|k\|^2)^{\frac{2-\varepsilon}{2}}\left(\frac{1}{\sqrt{-m'_k}}+\frac{1}{\|k\|}\right)^2\\
    &\le (1+\|k\|^2)^{\frac{2-\varepsilon}{2}}\left(\frac{1}{C(1+\|k\|^2)^{\frac{2-\varepsilon}{4}}}+\frac{1}{D(1+\|k\|^2)^{1/2}}\right)^2\\
    &= \frac{1}{C^2}+\frac{2}{CD(1+\|k\|^2)^{\varepsilon/4}}+\frac{1}{D^2(1+\|k\|^2)^{\varepsilon/2}}.
\end{align*}
This shows the uniform boundedness required to pass the limit inside the summation term-wise. The proof is then complete by using ~\eqref{col:beta}.
\end{proof}

\subsection{Spatial and temporal regularity over $L^2(\T^n)$}\label{sec:temporal1}
The solution $U(t)\in H^s(\T^n)$ given in \eqref{EqU-homo} is a distribution when $s<0$, and  when $s\ge 0$ defines a regular function 
\begin{align*}
    u(x,t):=U(t)(x)=(\hat{f}_0+\hat{g}_0t)+\sum_{0\ne k\in \mathbb{Z}^n}\left[\hat{f}_k\cos(\sqrt{-m_k}t)+\frac{\hat{g}_k}{\sqrt{-m_k}}\sin(\sqrt{-m_k}t)\right]e^{i k\cdot x}.
\end{align*}
Let the data be two regular functions; $f\in H^{s_1}(\T^n)$ and $g\in H^{s_2}(\T^n)$ with $s_1\ge 0$ and $s_2\ge 0$.
In this section, we focus on conditions that guarantee that $u(x,t)$, the solution of \eqref{eq:Homo-wave-general}, is a regular function $u(\cdot,t)\in L^2(\T^n)$. 
%\hl{sentence}
%This is also the unique solution to the  nonlocal wave equation \eqref{eq:Homo-wave-general} where the derivative with respect to time can be interpreted as Gateaux derivative with the $L^2$ norm or in the classical sense.

\subsubsection{Functions with absolutely summable Fourier coefficients}
%For completeness of the presentation, 
In this section, we state some well-know results on the differentiablity and summability of series in multidimensions. These results will be used in subsequent sections and are included here for completeness of the presentation.
%For $\zeta \in [0,1]$, a scalar function $f$ on $\mathbb{R}^n$ is called $\zeta$-H{\"o}lder continuous if there exists a constant $C$ such that 
%\begin{align*}
%    |f(x)-f(y)|\le C\|x-y\|^\zeta,
%\end{align*}
%for any $x,y$ in the domain of $f$.\\

%\noindent
The following theorem is a special case of Theorem 3.2.16 in  \cite{Grafakos2009}.
\begin{theorem}\label{Fourier_summable}
Suppose that $f$ is a function defined on the torus $\T^n$ satisfying the $\zeta$-H{\"o}lder continuity condition for some $\zeta>n/2$. Then, the Fourier coefficients of $f$ are absolutely summable.
\end{theorem}

%\subsection{Differentiation of a series}
The following fact on differentiating a series term by term will be used in subsequent sections. 
\begin{lemma}\label{diff}
Let $T$ be an open subset of $\mathbb{R}^+$, $K$ is a measurable space equipped with the counting measure and $f_k:T\to \mathbb{R}$. Suppose that
\begin{enumerate}
    \item $f_k(t)$ is summable over $k$ for each $t\in T$.
    \item For almost all $k\in K$, the derivative $\frac{d}{dt}f_k(t)$ exists and is continuous for all $t\in T$.
    \item There is a summable sequence $\theta_k$ such that $|\frac{d}{dt}f_k(t)|\le \theta_k$ for all $t\in T$ and almost all $k\in K$.
\end{enumerate}
Then, for all $t\in T$,
\begin{align*}
    \frac{d}{dt}\sum_{k}f_k(t)=\sum_{k}\frac{d}{dt}f_k(t).
\end{align*}
\end{lemma}

In addition, we include the following useful fact about the summability of certain series in multi-dimensions.
%\subsection{Some extra lemmas}
\begin{lemma}
\label{series_summability}
In $\mathbb{R}^n$, the series
\begin{equation}\label{series_convergence}
    \sum_{0\ne k\in \mathbb{Z}^n}\frac{1}{\|k\|^r}
\end{equation}
converges if and only if $r>n$.
\end{lemma}

\subsubsection{Temporal regularity with respect to the Gateaux derivative}
Let the data $f\in H^{s_1}(\T^n)$ and $g\in H^{s_2}(\T^n)$, with $s_1\ge 0$ and $s_2\ge 0$. 
\begin{theorem}
Let $n\ge 1$, $\delta>0$, and $\beta<n+4$. 
Suppose that
\begin{enumerate}
    \item $s_1,s_2\ge 0$, in the case when $\beta<n$,
    \item $s_1,s_2>0$, in the case when $\beta=n$, and
    \item $s_1\ge \frac{3}{2}(\beta-n) $ and $s_2\ge \beta-n$, in the case when $\beta>n$.
\end{enumerate}
Then, the  function $u$, where $  u(x,t)=U(t)(x)$, with $U$ defined  in \eqref{EqU-homo}, is the unique solution to the  nonlocal wave equation \eqref{eq:Homo-wave-general}.
Moreover, when $\beta\le n$, then $u(x,t)\in C^\infty([0,\infty); H^s(T^n))$, where $s=\min\{s_1,s_2\}$, and when $\beta>n$, 
%$u\in C^2([0,\infty); H^s(T^n))$, where
%s=\min\{s_1,s_2+\frac{\beta-n}{2}\}$.
then $u(x,t)\in C^2([0,\infty); H^s(T^n))$, where
$s=\min\{s_1,s_2+\frac{\beta-n}{2}\}$.
\end{theorem}

\subsubsection{Temporal regularity with respect to the classical derivative}
\begin{theorem}
Let $n\ge 1$, $\delta>0$ and $\beta<n$. Suppose that $f$ and $g$ satisfy the $\alpha$-Holder continuity condition with $\alpha>n/2$. Then $u(x,\cdot)\in C^\infty[0,\infty)$, for any $x\in \mathbf{T}^n$.
\end{theorem}
\begin{proof}
Fix a value of $p$ in $\mathbb{Z}_{\ge 0}$. Consider the two series $\sum\limits_{0\ne k\in\mathbb{Z}^n}\hat{f}_k\cos(\sqrt{-m_k}t)e^{i k\cdot x}$ and \\$\sum\limits_{0\ne k\in \mathbb{Z}^n}\frac{\hat{g}_k}{\sqrt{-m_k}}\sin(\sqrt{-m_k}t)e^{i k\cdot x}$. From Theorem~\eqref{thm:asymptotics}, $\sqrt{-m_k}$ is bounded for large value of $\|k\|$. Hence,
\begin{align*}
    &\left|\hat{f}_k\cos(\sqrt{-m_k}t)e^{i k\cdot x}\right|\le|\hat{f}_k|,\\
    &\left|\frac{d^p}{d t^p}(\hat{f}_k\cos(\sqrt{-m_k}t)e^{i k\cdot x})\right|\le (\sqrt{-m_k})^p|\hat{f}_k|\le C_1|\hat{f}_k|,\\
    &\left|\frac{\hat{g}_k}{\sqrt{-m_k}}\sin(\sqrt{-m_k}t)e^{i k\cdot x}\right|\le\frac{|\hat{g}_k|}{\sqrt{-m_k}}\le C_2|\hat{g}_k|,\\
    &\left|\frac{d^p}{d t^p}\left(\frac{\hat{g}_k}{\sqrt{-m_k}}\sin(\sqrt{-m_k}t)e^{i k\cdot x}\right)\right|\le(\sqrt{-m_k})^{p-1} |\hat{g}_k|\le C_3|\hat{g}_k|,
\end{align*}
as $\|k\|\to \infty$.  From \eqref{Fourier_summable}, the Fourier coefficients  $\hat{f}_k$ and $\hat{g}_k$ are absolutely summable. Now applying Lemma~\eqref{diff} completes the proof.
\end{proof}

\begin{theorem}
Let $n\ge 1$, $\delta>0$ and $\beta=n$. Suppose $s_1>n$ and $s_2>n$. Then $u(x,\cdot)\in C^\infty[0,\infty)$.
\end{theorem}

\begin{proof}
Fix a value of $p$ in $\mathbb{Z}_{\ge 0}$ and choose $\varepsilon>0$ such that $s_1-p\varepsilon/2>n$ and $s_2-(p-1)\varepsilon/2>n$. From \eqref{thm:asymptotics} we have that for large value of $\|k\|$, there exists constant $D,E>0$ such that
\begin{align*}
    -m_k\le D\log\|k\|\le E\|k\|^\varepsilon.
\end{align*}
Consider the two series $\sum\limits_{0\ne k\in\mathbb{Z}^n}\hat{f}_k\cos(\sqrt{-m_k}t)e^{i k\cdot x}$ and 
$\sum\limits_{0\ne k\in \mathbb{Z}^n}\frac{\hat{g}_k}{\sqrt{-m_k}}\sin(\sqrt{-m_k}t)e^{i k\cdot x}$. Note that from \eqref{thm:asymptotics}, $-m_k$ is bounded from below by some constant for large $\|k\|$  then use \eqref{Forier-decay} 
\begin{align*}
    &\left|\hat{f}_k\cos(\sqrt{-m_k}t)e^{i k\cdot x}|\le|\hat{f}_k\right|\le \frac{C_1}{\|k\|^{s_1}},\\
    &\left|\frac{d^p}{d t^p}(\hat{f}_k\cos(\sqrt{-m_k}t)e^{i k\cdot x})\right|\le (\sqrt{-m_k})^p|\hat{f}_k|\le\frac{C_2}{\|k\|^{s_1-p\varepsilon/2}},\\
    &\left|\frac{\hat{g}_k}{\sqrt{-m_k}}\sin(\sqrt{-m_k}t)e^{i k\cdot x}\right|\le\frac{|\hat{g}_k|}{\sqrt{-m_k}}\le\frac{C_3}{\|k\|^{s_2}\sqrt{-m_k}},\\
    &\left|\frac{d^p}{d t^p}\left(\frac{\hat{g}_k}{\sqrt{-m_k}}\sin(\sqrt{-m_k}t)e^{i k\cdot x}\right)\right|\le (\sqrt{-m_k})^{p-1}|\hat{g}_k|\le\frac{C_4}{\|k\|^{s_2-(p-1)\varepsilon/2}}.
\end{align*}
as $\|k\|\to\infty$. Applying Lemma \eqref{series_summability} and Lemma \eqref{diff} completes the proof.
\end{proof}

\begin{theorem}
Let $n\ge 1$, $\delta>0$ and $\beta>n$. Suppose that there exists $p\in\mathbb{Z}_{\ge 0}$ such that $s_1>n+p\frac{\beta-n}{2}$ and $s_2>n+(p-1)\frac{\beta-n}{2}$. Then, $u(x,\cdot)\in 
C^p[0,\infty)$, for any $x\in \mathbf{T}^n$.
\end{theorem}
\begin{proof}
From \eqref{thm:asymptotics}, for large $\|k\|$, there exists $D_1,D_2>0$ such that $D_1\|k\|^{\beta-n}\le-m_k\le D_2\|k\|^{\beta-n}$. Consider the two series $\sum\limits_{0\ne k\in\mathbb{Z}^n}\hat{f}_k\cos(\sqrt{-m_k}t)e^{i k\cdot x}$ and $\sum\limits_{0\ne k\in \mathbb{Z}^n}\frac{\hat{g}_k}{\sqrt{-m_k}}\sin(\sqrt{-m_k}t)e^{i k\cdot x}$ then use \eqref{Forier-decay}, we get
\begin{align*}
    &\left|\hat{f}_k\cos(\sqrt{-m_k}t)e^{i k\cdot x}|\le|\hat{f}_k\right|\le \frac{C_1}{\|k\|^{s_1}},\\
    &\left|\frac{d^p}{d t^p}(\hat{f}_k\cos(\sqrt{-m_k}t)e^{i k\cdot x})\right|\le (\sqrt{-m_k})^p|\hat{f}_k|\le\frac{C_2}{\|k\|^{s_1-p(\beta-n)/2}},\\
    &\left|\frac{\hat{g}_k}{\sqrt{-m_k}}\sin(\sqrt{-m_k}t)e^{i k\cdot x}\right|\le\frac{|\hat{g}_k|}{\sqrt{-m_k}}\le\frac{C_3}{\|k\|^{s_2+(\beta-n)/2}},\\
    &\left|\frac{d^p}{d t^p}\left(\frac{\hat{g}_k}{\sqrt{-m_k}}\sin(\sqrt{-m_k}t)e^{i k\cdot x}\right)\right|\le (\sqrt{-m_k})^{p-1}|\hat{g}_k|\le\frac{C_4}{\|k\|^{s_2-(p-1)(\beta-n)/2}}.
\end{align*}
The summability of these series follows from Lemma \eqref{series_summability}. Applying Lemma \eqref{diff} completes the proof. 
\end{proof}

\section{Nonlocal wave equation with forcing term}\label{sec:forcing}
In this section, we focus on the following nonlocal wave equation with a  forcing term and zero initial data. 
\begin{align}\label{eq:Nonhomo-wave-general}
    \begin{cases}
    u_{tt}(x,t)&=\Ldel u(x,t)+b(x),\quad x\in\T^n, \; t>0,\\
    u(x,0)&=0,\\
    u_t(x,0)&=0.
    \end{cases}
\end{align}
In a similar manner to Section \ref{sec:Homogeneous} ,the results of this section hold for any nonlocal operator $\Ldel$ with $\beta<n+4$. As being  recalled from Sections \ref{sec:prelim} and \ref{sec:L-explicit} that $\beta<n$ associates with an integrable kernel, $n\le \beta<n+2$ associates with a singular kernel, $\beta=n+2$ associates with $\Ldel=\Delta$ (the Laplacian), and $n+2<\beta<n+4$ associates with an integro-differential operator as given in \eqref{eq:nonlocal_laplacian_extended}.  \\
In order to study the existence, uniqueness, and regularity of solutions to \eqref{eq:Nonhomo-wave-general} over the space of periodic distributions, we consider the identification $U(t)=u(\cdot,t)$, with $U:[0,\infty)\rightarrow H^s(\T^n)$.
\subsection{Spatial regularity of the solution}
Fix $b\in H^{\sigma}(\T^n)$ for $\sigma\in\mathbb{R}$. For any $t\ge 0$, define
\begin{align}\label{EqU-forcing}
    U(t):=\sum_{k\in \mathbb{Z}^n}\hat{U}_k(t)e^{ik\cdot x}=\hat{b}_0\frac{t^2}{2}+\sum_{0\ne k\in \mathbb{Z}^n}\frac{\hat{b}_k}{m_k}\left[\cos(\sqrt{-m_k}t)-1\right]e^{i k\cdot x}.
\end{align}
\begin{theorem}\label{thm:spatial_reg2}
Let $n\ge 1$, $\delta>0$, and $\beta<n+4$. For any $t\ge 0$, $U(t)\in H^{s}(\T^n)$ where \\$s=\sigma+\max\{0,\beta-n\}$.
\end{theorem}
\begin{proof}
Denote $\theta=\max\{0,\beta-n\}$. Observe that
\begin{align*}
    \sum_{k\in\mathbb{Z}^n}(1+\|k\|^2)^{s}|\hat{U}_k(t)|^2&=|\hat{b}_0|\frac{t^2}{2}+ \sum_{0\ne k\in \mathbb{Z}^n}(1+\|k\|^2)^{s}\frac{{\hat{b}_k}^2}{m^2_k}\left[\cos(\sqrt{-m_k}t)-1\right]^2\\
    &=|\hat{b}_0|\frac{t^2}{2}+ \sum_{0\ne k\in \mathbb{Z}^n}(1+\|k\|^2)^{\theta}\frac{\left[\cos(\sqrt{-m_k}t)-1\right]^2}{m^2_k}(1+\|k\|^2)^{\sigma}{\hat{b}_k}^2.
\end{align*}
Since $b\in H^\sigma$, the proof will be complete provided that
\begin{align*}
    \frac{(1+\|k\|^2)^{\theta}}{m^2_k},
\end{align*}
is bounded for large values of $\|k\|$.\\
To prove this fact, we consider two cases. First, when $\beta\le n$, then $\theta=0$ and by Theorem ~\eqref{thm:asymptotics} there exist $r_1>0$ and $C_1>0$ such that   $|m_k|\ge C_1$ for $\|k\|\ge r_1$. This leads us to
\begin{align*}
    \frac{(1+\|k\|^2)^{\theta}}{m^2_k}\le \frac{1}{C^2_1}.
\end{align*}
With the same approach, for $\beta> n$, then $\theta=\beta-n$ and by Theorem \eqref{thm:asymptotics}  there exist $r_2>0$ and $C_2>0$ such that  $|m_k|\ge C_2\|k\|^{\beta-n}$, for $\|k\|\ge r_2$. This leads  to
\begin{align*}
    \frac{(1+\|k\|^2)^{\theta}}{m^2_k}\le \frac{1}{C_2^2}\left(\frac{1+\|k\|^2}{\|k\|^2}\right)^{\beta-n},
\end{align*}
which is bounded.
\end{proof}

\begin{theorem}\label{thm:regularity1}
Let $n\ge 1$, $\delta>0$, and $\beta<n+4$. Let $U(t)$ be the map given by \eqref{EqU-forcing}. Thus,
\begin{enumerate}
    \item if $\beta<n$, then $U(t)\in C^\infty([0,\infty),H^q(\T^n))$ for any $q\le \sigma$,
    \item if $\beta=n$, then $U(t)\in C^\infty([0,\infty),H^q(\T^n))$ for any $q<\sigma$, and
    \item if $\beta>n $, then $U(t)\in C^{p+1}([0,\infty),H^q(\T^n))$ for any positive integer $p$ satisfying 
    \begin{align*}
        q-\sigma+(p-1)\frac{\beta-n}{2}\le 0.    \end{align*}
    
\end{enumerate}
\end{theorem}

\begin{proof}
We will show this result using induction. Suppose $U$ is already differentiable up to $p$ times, then by Lemma \eqref{lem:derivative}, we should have
\begin{align*}
    U^{(p)}(t)=\sum_{k\in\mathbb{Z}^n}\hat{U}^{(p)}_k(t)e^{ik\cdot x},
\end{align*}
where $\hat{U}^{(p)}_k(t)$ is given by
\begin{align*}
    \hat{U}^{(p)}_k(t)=-\hat{b}_k\sqrt{-m_k}^{p-2}\cos\left(\sqrt{-m_k}t+\frac{p\pi}{2}\right),
\end{align*}
for $k\ne 0$ and $\hat{U}_0^{(1)}=\hat{b}_0 t$,\;$\hat{U}_0^{(2)}=\hat{b}_0$ and $\hat{U}_0^{(p)}=0$ for all $p\ge 3$. We will show that 
\begin{align*}
    U^{(p+1)}(t)=\sum_{k\in\mathbb{Z}^n}\hat{U}^{(p+1)}_k(t)e^{ik\cdot x}
    .
\end{align*}
Observe that using the mean value theorem
\begin{align*}
    T_h:&=\left\|\frac{U^{(p)}(t+h)-U^{(p)}(t)}{h}-\sum_{k\in\mathbb{Z}^n}\hat{U}^{(p+1)}_k(t)e^{ik\cdot x}\right\|^2_{H^q(\T^n)}\\
    &=\sum_{k\in\mathbb{Z}^n}(1+\|k\|^2)^q\left|\frac{\hat{U}^{(p)}_k(t+h)-\hat{U}^{(p)}_k(t)}{h}-\hat{U}^{p+1}_k(t)\right|^2,\\
    &=\left|\frac{\hat{U}^{(p)}_0(t+h)-\hat{U}^{(p)}_0(t)}{h}-\hat{U}^{p+1}_0(t)\right|^2\\
    &\qquad+\sum_{0\ne k\in\mathbb{Z}^n}(1+\|k\|^2)^q|\hat{b}_k|^2(-m_k)^{p-1}\left|\cos\left(\sqrt{-m_k}(t+\xi_p)+\frac{(p+1)\pi}{2}\right)\right.\\
    &\qquad\qquad\qquad\left.-\cos\left(\sqrt{-m_k}t+\frac{(p+1)\pi}{2}\right)\right|^2,
\end{align*}
for some $\xi_{p,k}\in(0,h)$. This means that as $h\to 0$, each term of the summation will approach 0. The proof will be complete if we are able to pass the limit inside the summation. This is possible provided the uniform boundedness of
\begin{align*}
    a_k:=(1+\|k\|^2)^{q-\sigma}(-m_k)^{p-1}.
\end{align*}
To establish the bound for these terms, we will consider three cases when $\beta<n$, $\beta=n$ and $\beta>n$,  respectively.\\
When $\beta<n$, from Theorem \eqref{thm:asymptotics}, $-m_k$ is bounded, and this  provides the uniform boundedness of $a_k$ provided that $q\le \sigma$.\\
When $\beta=n$, fix an $\epsilon>0$. From Theorem \eqref{thm:asymptotics}, there exists $C,r>0$ such that $-m_k\le C\log\|k\|\le C(1+\|k\|^2)^{\epsilon/2}$ whenever $\|k\|\ge r$. Then for large $\|k\|$,
\begin{align*}
    a_k\le C(1+\|k\|^2)^{q-\sigma+(p-1)\epsilon/2}.
\end{align*}
This yields the uniform boundedness of $a_k$ as long as $q-\sigma+(p-1)\epsilon/2\le 0$. So given $q<\sigma$, we can always pick appropriate $\epsilon$ to satisfy this.\\
When $\beta>n$, from Theorem \eqref{thm:asymptotics}, there is $C,r>0$ such that $-m_k\le C\|k\|^{\beta-n}\le C(1+\|k\|^2)^{(\beta-n)/2}$ for $\|k\|\ge r$. Then for large $\|k\|$
\begin{align*}
    a_k\le C(1+\|k\|^2)^{q-\sigma+(p-1)\frac{\beta-n}{2}}.
\end{align*}
This shows that the uniform boundedness of $a_k$ is established given the conditions
\begin{align*}
    q-\sigma+(p-1)\frac{\beta-n}{2}\le 0.
\end{align*}
\end{proof}

 The following theorem is the direct application of Theorem \eqref{thm:spatial_reg2} and Theorem \eqref{thm:regularity1}. The assumptions are to guarantee that the second derivative exists.
\begin{theorem}
Let $n\ge1$, $\delta>0$ and $\beta<n+4$. Let $b\in H^{\sigma}(\T^n)$ for $\sigma\in\mathbb{R}$ and suppose that
\begin{enumerate}
    \item $q\le \sigma$ in the case when $\beta<n$,
    \item $q<\sigma$ in the case when $\beta=n$, and
    \item $q\le \sigma$ in the case when $\beta>n$.
\end{enumerate}
Then, the map $U(t)$ defined by \eqref{EqU-forcing} is the unique solution in the space $H^q(\T^n)$ of the nonlocal wave equation
\begin{align}\label{eq:Forcing-wave}
    \begin{cases}
    \frac{d^2}{dt^2}U(t)&=\Ldel U(t)+b,\\
    U(0)&=0,\\
    \frac{d}{dt}U(0)&=0.
    \end{cases}
\end{align}
Moreover, when $\beta\le n$, $U(t)\in C^\infty([0,\infty); H^{\sigma}(\T^n))$  and when $\beta>n$, 
$U(t)\in C^2([0,\infty); H^{\sigma+\beta-n}(\T^n))$.
\end{theorem}

\subsection{Convergence to the classical wave equation solution}
\label{sec:convergence2}
In this section, we will provide some results on the convergence of the solutions of the nonlocal wave equation \eqref{eq:Forcing-wave} to the local solution for two types of  limits: as $\delta \to 0^+$, with $\beta< n+4$ is being fixed,  or as $\beta\to n+2$, with $\delta>0$ is being fixed.\\
The local wave equation is defines as
\begin{align}\label{eq:Forcing-wave-local}
    \begin{cases}
    \frac{d^2}{dt^2}U^0(t)&=\Delta U^0(t)+b,\\
    U^0(0)&=0,\\
    \frac{d}{dt}U^0(0)&=0,
    \end{cases}
\end{align}
where its solution is given by
\begin{align*}
    U^0(t):=\hat{b}_0\frac{t^2}{2}+\sum_{0\ne k\in \mathbb{Z}^n}\frac{\hat{b}_k}{-\|k\|^2}\left[\cos(\|k\|t)-1\right]e^{i k\cdot x}.
\end{align*}

To emphasize the dependence of the solution of the nonlocal wave equation \eqref{eq:Homo-wave} on the parameters $\delta$ and $\beta$, in the remaining part of this section we denote the solution $U$ by $U^{\delta,\beta}$.
\begin{theorem}
Let $n\ge 1$, $\delta>0$ and $\beta<n+4$. Suppose $b\in H^{\sigma}(\T^n)$ for $\sigma\in\mathbb{R}$. Let $U^0(t)\in H^{s'}(\T^n)$ and $U^{\delta,\beta}(t)\in H^{s}(\T^n)$ be the solutions to the local and nonlocal wave equation 
 \eqref{eq:Forcing-wave-local} and \eqref{eq:Forcing-wave} respectively, where $s'=\sigma+2$ and $s=\sigma+\max\{0,\beta-n\}$. Then when $\beta \le n+2$
\begin{equation*}
    \lim_{\delta\to 0^+}U^{\delta,\beta}(t)=U(t) \quad \text{in} \quad H^{s}(\T^n),
\end{equation*}
and when $n+2<\beta<n+4$
\begin{equation*}
    \lim_{\delta\to 0^+}U^{\delta,\beta}(t)=U(t) \quad \text{in} \quad H^{s'}(\T^n).
\end{equation*}
\end{theorem}

\begin{proof} 
It is easy to check that $H^{s'}(\T^n)\subset H^{s}(\T^n)$ when $\beta \le n+2$ and $H^{s}(\T^n)\subset H^{s'}(\T^n)$ when $n+2<\beta< n+4$. Hence all the limits make sense. From \eqref{EqU-forcing}, the Fourier coefficients $\hat{U}_k^{\delta,\beta}(t)$ in the nonlocal case and $\hat{U}_k(t)$ in the classical case ($\beta=n+2$) are determined by
\begin{align*}
    \hat{U}_k^{\delta,\beta}(t)&=\frac{\hat{b}_k}{m_k}\left[\cos(\sqrt{-m_k}t)-1\right],\\
    \hat{U}_k(t)&=\frac{\hat{b}_k}{-\|k\|^2}\left[\cos(\|k\|t)-1\right],
\end{align*}
for nonzero $k$ and $\hat{U}_0^{\delta,\beta}(t)=\hat{U}_0(t)=\hat{b}_0t^2/2$. From these observation, for $s_0=s$ when $\beta\le n+2$ and $s_0=s'$ when $n+2<\beta<n+4$, one can get
\begin{align*}
    \|U^{\delta,\beta}(t)-U(t) \|^2_{H^{s_0}(\T^n)} &= \sum_{0\ne k\in\mathbb{Z}^n}(1+\|k\|^2)^{s_0}\|\hat{U}^{\delta,\beta}_k(t) - \hat{U}_k(t)\|^2 \\
&=\sum_{0\ne k\in\mathbb{Z}^n}(1+\|k\|^2)^{s_0}|\hat{b}_k|^2\left|\frac{\cos(\sqrt{-m_k}t)-1}{m_k}-\frac{\cos(\sqrt{\|k\|}t)-1}{-\|k\|^2}\right|^2.
\end{align*}
Since $b\in H^{\sigma}(\T^n)$, one can pass the limit as $\delta\to 0^+$ term by term inside the summation provided the uniform boundedness of
\begin{align*}
    (1+\|k\|^2)^{s_0-\sigma}\left|\frac{\cos(\sqrt{-m_k}t)-1}{m_k}-\frac{\cos(\sqrt{\|k\|}t)-1}{-\|k\|^2}\right|^2,
\end{align*}
for large $\|k\|$. Denote $\theta=\max\{0,\beta-n\}$. Since $s_0-\sigma=\theta$ when $\beta\le n+2$ and \\$s_0-\sigma=2\le \theta$ when $n+2<\beta<n+4$, one obtains
\begin{align*}
    (1+\|k\|^2)^{s_0-\sigma}\left|\frac{\cos(\sqrt{-m_k}t)-1}{m_k}-\frac{\cos(\sqrt{\|k\|}t)-1}{-\|k\|^2}\right|^2\\
    \le 4(1+\|k\|^2)^{s_0-\sigma}\left(\frac{1}{-m_k}+\frac{1}{\|k\|^2}\right)^2.
\end{align*}
We consider two cases. When $\beta\le n$, then $s_0-\sigma=0$ . From Theorem \eqref{thm:asymptotics}, there exists $r_1>0$ and $C_1>0$ such that $-m_k\ge C_1$ for $\|k\|\ge r_1$. When $\beta>n$, then $s_0-\sigma\le\beta-n$. From ~\eqref{thm:asymptotics}, there exists $r_2>0$ and $C_2>0$ such that $-m_k\ge C_2\|k\|^{\beta-n}$ for $\|k\|\ge r_2$. These results lead us to the uniform boundedness of the right hand side.\\
Now, taking the limit as $\delta\to0^+$ inside the summation term-wise and applying the pointwise convergence of ~\eqref{col:delta} completes the proof.
\end{proof}

\begin{theorem}\label{converge2}
Let $n\ge 1$, $\delta>0$. Suppose $b\in H^{\sigma}(\T^n)$ for $\sigma\in\mathbb{R}$. Let $U(t)\in H^{s'}(\T^n)$ and $U^{\delta,\beta}(t)\in H^{s}(\T^n)$ be the solution to the local and nonlocal wave equation respectively, where $s'=\sigma+2$ and $s=\sigma+\max\{0,\beta-n\}$. Then for any $\varepsilon\in(0,2)$, $\varepsilon<n/2$
\begin{equation*}
    \lim_{\beta\to n+2}U^{\delta,\beta}(t)=U(t) \quad \text{in} \quad H^{s_0}(\T^n),
\end{equation*}
where $s_0=\sigma+2-\varepsilon$.
\end{theorem}
\begin{proof}
Fix an $\varepsilon \in (0,2),\varepsilon<n/2$ and define $\beta'=n+2-\varepsilon>\frac{n+4}{2}$. 
For simplicity, we would still denote 
\begin{align*}
    m_k:=m^{\delta,\beta}(k)\quad\text{and}\quad m'_k:=m^{\delta,\beta'}(k).
\end{align*}
Then, for any $\beta$ such that $|\beta-(n+2)|<\varepsilon$ or $\frac{n+4}{2}<\beta'<\beta<n+2+\varepsilon<n+4$, we have  $\frac{\beta-n}{2}>\frac{2-\varepsilon}{2}$. This shows that $s,s'>s_0$, so the limit in the statement of the theorem makes sense.
Now consider
\begin{align*}
    \|U^{\delta,\beta}(t)-U(t) \|^2_{H^{s_0}(\T^n)} &=\sum_{0\ne k\in\mathbb{Z}^n}(1+\|k\|^2)^{s_0}|\hat{b}_k|^2\left|\frac{\cos(\sqrt{-m_k}t)-1}{m_k}-\frac{\cos(\sqrt{\|k\|}t)-1}{-\|k\|^2}\right|^2.
\end{align*}
Similar to the previous theorem, we need to establish the uniform boundedness of
\begin{align*}
    (1+\|k\|^2)^{s_0-\sigma}\left|\frac{\cos(\sqrt{-m_k}t)-1}{m_k}-\frac{\cos(\sqrt{\|k\|}t)-1}{-\|k\|^2}\right|^2,
\end{align*}
for large $\|k\|$. For this term, from ~\eqref{thm:asymptotics} there exists $r>0$ and $C,D>0$ such that\\ $-m'_k\ge C(1+\|k\|^2)^{\frac{2-\varepsilon}{2}}$ and $\|k\|^2\ge D(1+\|k\|^2)$ for $\|k\|\ge r$. This yields
\begin{align*}
    &(1+\|k\|^2)^{s_0-\sigma}\left|\frac{\cos(\sqrt{-m_k}t)-1}{m_k}-\frac{\cos(\sqrt{\|k\|}t)-1}{-\|k\|^2}\right|^2\\
    &\le
    4(1+\|k\|^2)^{s_0-\sigma}\left(\frac{1}{-m_k}+\frac{1}{\|k\|^2}\right)^2\\
    &\le 4(1+\|k\|^2)^{2-\varepsilon}\left(\frac{1}{-m'_k}+\frac{1}{\|k\|^2}\right)^2\\
    &\le (1+\|k\|^2)^{2-\varepsilon}\left(\frac{1}{C(1+\|k\|^2)^{\frac{2-\varepsilon}{2}}}+\frac{1}{D(1+\|k\|^2)}\right)^2\\
    &= \frac{1}{C^2}+\frac{2}{CD(1+\|k\|^2)^{\varepsilon/2}}+\frac{1}{D^2(1+\|k\|^2)^{\varepsilon}}.
\end{align*}
This shows the uniform boundedness required to pass the limit inside the summation term-wise. The proof is then complete by using ~\eqref{col:beta}.
\end{proof}

\subsection{Spatial and temporal regularity over $L^2(\T^n)$}\label{sec:temporal2}
The solution $U(t)\in H^s(\T^n)$ given in \eqref{EqU-forcing} is a distribution when $s<0$, and  when $s\ge 0$ defines a regular function
\begin{align*}
    u(x,t):=U(t)(x)=\hat{b}_0\frac{t^2}{2}+\sum_{0\ne k\in \mathbb{Z}^n}\frac{\hat{b}_k}{m_k}\left[\cos(\sqrt{-m_k}t)-1\right]e^{i k\cdot x}.
\end{align*}
Let the forcing term be a regular functions $b\in H^{\sigma}(\T^n)$  with $\sigma\ge 0$.
In this section, we focus on conditions that guarantee that $u(x,t)$, the solution of \eqref{eq:Forcing-wave}, is a regular function $u(\cdot,t)\in L^2(\T^n)$.

\subsubsection{Temporal regularity with respect to the Gateaux derivative} Let the forcing term $b\in H^\sigma(\T^n)$ with $\sigma\ge 0$.
\begin{theorem}
Let $n\ge 1$, $\delta>0$ and $\beta<n+4$. Suppose that
\begin{enumerate}
    \item $s\ge 0$ in the case when $\beta<n$,
    \item $s>0$ in the case when $\beta=n$, and
    \item $s\ge 0 $ in the case $\beta>n$.
\end{enumerate}
Then the  function $u$, where $  u(x,t)=U(t)(x)$, with $U$ defined  in \eqref{EqU-forcing}, satisfies the classical nonlocal wave equation \eqref{eq:Forcing-wave}.
Moreover, when $\beta\le n$, $u(x,t)\in C^\infty([0,\infty); H^\sigma(\T^n))$, and when $\beta>n$,
$u(x,t)\in C^2([0,\infty); H^{\sigma+\beta-n}(\T^n))$.
\end{theorem}

\subsubsection{Temporal regularity with respect to the classical derivative}
\begin{theorem}
Let $n\ge 1$, $\delta>0$ and $\beta<n$. Suppose that $b$ satisfy the $\gamma$-H{\"o}lder continuity condition with $\gamma>n/2$. Then $u(x,\cdot)\in C^\infty[0,\infty)$, for any $x\in \mathbf{T}^n$.
\end{theorem}
\begin{proof}
Fix a value of $p$ in $\mathbb{Z}_{\ge 0}$. Consider the series 
\[
\sum_{0\ne k\in\mathbb{Z}^n}\frac{\hat{b}_k}{m_k}\left[\cos(\sqrt{-m_k}t)-1\right]e^{i k\cdot x}.
\] 
From ~\eqref{thm:asymptotics}, $\sqrt{-m_k}$ is bounded for large value of $\|k\|$. Hence,
\begin{align*}
    &\left|\frac{\hat{b}_k}{m_k}\left[\cos(\sqrt{-m_k}t)-1\right]e^{i k\cdot x}\right|\le\frac{|\hat{b}_k|}{|m_k|}\le C_1|\hat{b}_k|,\\
    &\left|\frac{d^p}{d t^p}\frac{\hat{b}_k}{m_k}\left[\cos(\sqrt{-m_k}t)-1\right]e^{i k\cdot x}\right|\le (\sqrt{-m_k})^{p-2}|\hat{b}_k|\le C_2|\hat{b}_k|
\end{align*}
as $\|k\|\to \infty$.  From \eqref{Fourier_summable}, the Fourier coefficients  $\hat{b}_k$ is absolutely summable. Now applying Lemma \eqref{diff} completes the proof.
\end{proof}

\begin{theorem}
Let $n\ge 1$, $\delta>0$ and $\beta=n$. Suppose $\sigma>n$. Then $u(x,\cdot)\in C^\infty[0,\infty)$.
\end{theorem}

\begin{proof}
Fix a value of $p$ in $\mathbb{Z}_{\ge 0}$ and choose $\varepsilon>0$ such that $\sigma-\varepsilon(p-2)/2>n$. From ~\eqref{thm:asymptotics} we have that for large value of $\|k\|$, there exists constant $D,E>0$ such that
\begin{align*}
    -m_k\le D\log\|k\|\le E\|k\|^\varepsilon.
\end{align*}
Consider the series $\sum\limits_{0\ne k\in\mathbb{Z}^n}\frac{\hat{b}_k}{m_k}\left[\cos(\sqrt{-m_k}t)-1\right]e^{i k\cdot x}$. Note that from ~\eqref{thm:asymptotics}, $-m_k$ is bounded from below by some constant for large $\|k\|$  then use ~\eqref{Forier-decay} 
\begin{align*}
    &\left|\frac{\hat{b}_k}{m_k}\left[\cos(\sqrt{-m_k}t)-1\right]e^{i k\cdot x}\right|\le\frac{|\hat{b}_k|}{|m_k|}\le C_1|\hat{b}_k|,\\
    &\left|\frac{d^p}{d t^p}\frac{\hat{b}_k}{m_k}\left[\cos(\sqrt{-m_k}t)-1\right]e^{i k\cdot x}\right|\le (\sqrt{-m_k})^{p-2}|\hat{b}_k|\le \frac{C_2}{\|k\|^{\sigma-\varepsilon(p-2)/2}}
\end{align*}
as $\|k\|\to\infty$. The summability of these series follows from Lemma \eqref{series_summability}. Applying Lemma \eqref{diff} completes the proof. 
\end{proof}

\begin{theorem}
Let $n\ge 1$, $\delta>0$ and $\beta>n$. Suppose that there exists $p\in\mathbb{Z}_{\ge 0}$ such that $\sigma>n+(p-2)\frac{\beta-n}{2}$. Then $u(x,\cdot)\in C^p[0,\infty)$, for any $x\in \mathbf{T}^n$.
\end{theorem}
\begin{proof}
From \eqref{thm:asymptotics}, for large $\|k\|$, there exists $D_1,D_2>0$ such that $D_1\|k\|^{\beta-n}\le-m_k\le D_2\|k\|^{\beta-n}$. Consider the series $\sum\limits_{0\ne k\in\mathbb{Z}^n}\frac{\hat{b}_k}{m_k}\left[\cos(\sqrt{-m_k}t)-1\right]e^{i k\cdot x}$ then use \eqref{Forier-decay}, we get
\begin{align*}
    &\left|\frac{\hat{b}_k}{m_k}\left[\cos(\sqrt{-m_k}t)-1\right]e^{i k\cdot x}\right|\le\frac{|\hat{b}_k|}{|m_k|}\le \frac{C_1}{\|k\|^{\sigma+(\beta-n)}},\\
    &\left|\frac{d^p}{d t^p}\frac{\hat{b}_k}{m_k}\left[\cos(\sqrt{-m_k}t)-1\right]e^{i k\cdot x}\right|\le (\sqrt{-m_k})^{p-2}|\hat{b}_k|\le \frac{C_2}{\|k\|^{\sigma-(p-2)(\beta-n)/2}}
\end{align*}
Applying Lemma \eqref{series_summability} and Lemma \eqref{diff} completes the proof. 
\end{proof}

\begin{comment}
\section{Conclusion}
\begin{theorem}
Let $b\in H^{S}(\T^n)$,$f\in H^{s_1}(\T^n)$ and $g\in H^{s_2}(\T^n)$ such that $\hat{f}_0=0$. Suppose that
\begin{enumerate}
    \item $q\le s_1,s_2,S$ in the case $\beta<n$,
    \item $q<s_1,s_2,S$ in the case $\beta=n$,
    \item $q\le S$, $q-s_1+(\beta-n)\le 0$ and $q-s_2+(\beta-n)/2\le 0$ in the case $\beta>n$.
\end{enumerate}
Then the map $U(t)$ defined by
\begin{align*}
     U(t)=\hat{g}_0t+\hat{b}_0\frac{t^2}{2}+\sum_{0\ne k\in \mathbb{Z}^n}\left[\hat{f}_k\cos(\sqrt{-m_k}t)+\frac{\hat{g}_k}{\sqrt{-m_k}}\sin(\sqrt{-m_k}t)+\frac{\hat{b}_k}{m_k}\left(\cos(\sqrt{-m_k}t)-1\right)\right]e^{i k\cdot x}.
\end{align*}
is the unique solution in the space $H^q(\T^n)$ of the Nonlocal Wave equation
\begin{align*}
    \frac{d^2}{dt^2}U(t)&=\Ldel U(t)+b,\\
    U(0)&=f,\\
    \frac{d}{dt}U(0)&=g.
\end{align*}
Moreover, when $\beta\le n$, $U(t)\in C^\infty([0,\infty); H^s(\T^n))$ where $s=\min\{s_1,s_2,S\}$,  and when $\beta>n$, 
$U(t)\in C^2([0,\infty); H^s(\T^n))$ where $s=\min\{s_1,s_2+\frac{\beta-n}{2},S+\beta-n\}$.

\end{theorem}

\end{comment}

\clearpage
\bibliographystyle{acm}
\bibliography{refs}

\end{document}